\documentclass[preprint,12pt,longtitle]{elsarticle}

\usepackage[T1]{fontenc}
\usepackage[utf8]{inputenc}
\DeclareUnicodeCharacter{FF01}{!}
\usepackage{amsmath,amssymb,mathtools,mathrsfs}
\usepackage{graphicx}
\usepackage{booktabs,longtable,array,calc,tabularx}
\usepackage[section]{placeins}
\usepackage{flafter}
\usepackage{float}
\usepackage{needspace}
\usepackage[ruled,linesnumbered]{algorithm2e}
\usepackage{microtype}
\usepackage{xurl}
\usepackage[hidelinks]{hyperref}

\biboptions{sort&compress}
\graphicspath{{./}}
\allowdisplaybreaks[2]
\newcommand{\TableFont}{\fontsize{10}{12}\selectfont}
\AtBeginEnvironment{longtable}{\TableFont\setlength{\tabcolsep}{4pt}\renewcommand{\arraystretch}{1.12}}
\AtBeginEnvironment{table}{\TableFont}
\renewcommand{\topfraction}{0.90}

\renewcommand{\textfraction}{0.10}
\renewcommand{\floatpagefraction}{0.82}
\SetAlFnt{\fontsize{10}{12}\selectfont}
\SetAlCapFnt{\small}
\SetAlCapNameFnt{\small}
\SetAlgoSkip{medskip}

\SetKwInput{KwIn}{Input}
\SetKwInput{KwOut}{Output}

\begin{document}

\begin{frontmatter}

\title{Quantum Frozen--Oseen homotopy analysis method with LCHS for solving nonlinear partial differential equations}

\author[aff1]{Jinhao Yuan}
\author[aff1]{Leyu Chen}
\author[aff1,aff2]{Tiegang Liu}
\author[aff1,aff2]{Kun Wang\corref{cor1}}
\cortext[cor1]{Corresponding author}
\ead{wangkun@buaa.edu.cn}

\affiliation[aff1]{
  organization={LMIB and School of Mathematical Sciences},
  addressline={Beihang University},
  city={Beijing},
  postcode={100191},
  country={P.R. China}
}

\affiliation[aff2]{
  organization={International Research Center for Mathematics and Interdisciplinary Sciences},
  addressline={Hangzhou International Innovation Institute of Beihang University},
  city={Hangzhou},
  postcode={311115},
  country={P.R. China}
}

\begin{abstract}
Nonlinear partial differential equations (PDEs) underpin computational fluid dynamics, yet resolving nonlinear transport on fine grids remains computationally demanding. Quantum linear-evolution algorithms offer a possible route to large-scale simulation but cannot directly propagate nonlinear coupling. Here we develop a Frozen--Oseen quantum homotopy method (FOQHAM) framework that links transport-aware auxiliary-operator selection to the size of the resulting linear representation. We freeze the full Fréchet derivative at a prescribed flow profile, retaining transport and profile-gradient coupling and close each prescribed finite-order homotopy hierarchy as a fixed affine linear system by product lifting. We formulate its propagation through Linear Combinations of Hamiltonian Simulations (LCHS), without outer homotopy iterations or profile updates. For spatially semidiscrete equations, we establish local approximation bounds and exact closure of the finite-order hierarchy. Classical tests on Burgers, Korteweg--de Vries, Zakharov--Kuznetsov, and one- and two-dimensional isothermal compressible Navier--Stokes equations provide numerical evidence of accurate non-iterative approximations in the tested regimes. The framework thus connects nonlinear PDE approximation to quantum linear evolution and offers a conditional path toward quantum acceleration.
\end{abstract}


\end{frontmatter}

\section{Introduction}\label{introduction}

Computational fluid dynamics (CFD) connects fluid physics to predictions used in aerodynamic and multiphysics design \cite{ref64,ref65}. Its numerical core is the solution of nonlinear partial differential equations (PDEs), including conservation laws and the Navier--Stokes equations \cite{ref38,ref64,ref65}. Thin layers and interacting scales demand fine grids, while stiffness and long integration intervals increase computational and storage costs \cite{ref44,ref49,ref56,ref57,ref58}. These demands motivate alternative representations of large discretized systems.

Quantum computing encodes selected high-dimensional states in amplitudes and provides algorithms for linear systems and differential equations \cite{ref1,ref3,ref4,ref28,ref29}. Quantum spectral and finite-element methods connect these algorithms to PDE discretizations \cite{ref2,ref30,ref31}. Hamiltonian simulation implements unitary evolution \cite{ref32,ref33}; Schrödingerization embeds nonunitary linear dynamics in Hamiltonian systems \cite{ref66}. Linear Combinations of Hamiltonian Simulations (LCHS) instead represent nonunitary propagators as weighted unitary evolutions \cite{ref15,ref16}. Nonlinear transport, however, couples unknown components through products that these linear primitives cannot propagate directly \cite{ref5,ref6,ref7}.

Carleman methods address this difficulty by lifting monomials into an infinite linear hierarchy, with problem-dependent truncation bounds \cite{ref5,ref6,ref8,ref34}. Liouville and Koopman formulations evolve distributions or observables, but finite-dimensional closure is problem dependent \cite{ref7,ref35}. The homotopy analysis method (HAM) constructs a deformation series using an auxiliary linear operator and a convergence-control parameter \cite{ref9,ref10,ref11,ref36,ref37}. Xue et al.'s quantum HAM (QHAM) lifts retained coefficient products into tensor variables, closing a prescribed finite hierarchy as a linear system \cite{ref12}. Transport omitted from the auxiliary operator must be recovered through higher-order corrections or outer iteration. The resulting question is how auxiliary selection balances approximation order, representation size and implementation cost.

Here we examine this balance through a Frozen--Oseen homotopy construction. Classical Oseen and Picard linearizations incorporate prescribed flow information into linear transport operators \cite{ref13,ref14}. We freeze the full Fr\'echet derivative at a prescribed profile, retaining transport and profile-gradient coupling together with the residual and nonlinear remainder. This choice relates to exponential Rosenbrock linearization \cite{ref56}, but keeps the profile and lifted generator fixed. For quadratic equations, it specializes Xue et al.'s QHAM while retaining its product-lifting principle. Local approximation bounds and exact finite-order closure connect sufficient homotopy order to representation size for spatially semidiscrete polynomial systems. Through our construction, we use LCHS to propagate the system without outer homotopy iteration and profile updates.

Numerical tests cover Burgers, Korteweg--de Vries, Zakharov--Kuznetsov, and one- and two-dimensional isothermal compressible Navier--Stokes systems. Separate coefficient controls hold discretization and initialization fixed to isolate auxiliary selection from profile enrichment and other recipe choices. They show that more retained linearized physics can lower homotopy order, but need not minimize representation size. Spatial and finite-rule diagnostics distinguish approximation from propagation error; resource analysis retains matrix-access, normalization, source-integration and readout costs. Sections 2--4 develop the method, resource conditions and numerical evidence; Sections 5--6 discuss its scope and conclude.

\section{Quantum Frozen--Oseen homotopy analysis method}\label{quantum-frozen-oseen-homotopy-analysis-method}

\subsection{Frozen--Oseen homotopy analysis}\label{frozen-oseen-homotopy-analysis}

\subsubsection{Problem setting and boundary treatment}\label{problem-setting-and-boundary-treatment}

On a real or complex Hilbert space \(X\) with homogenized boundaries, consider the autonomous problem

\[
\frac{\mathrm d u}{\mathrm dt}
=\mathcal F(u)+f,\qquad
u(0)=u_{\mathrm{in}},\qquad
0\le t\le T,
\tag{1}
\]

Here \(\mathcal F:\mathcal D(\mathcal F)\subset X\to X\), and \(f\) is time-independent.

For a polynomial nonlinearity of maximum degree \(d\),

\[
\mathcal F(u)=\mathcal A u+
\sum_{r=2}^{d}\mathcal B_r(u,\ldots,u),
\tag{2}
\]

where \(\mathcal A\) is a closed linear differential operator and \(\mathcal B_r\) is an \(r\)-linear differential product on a common regularity domain. For example,

\[
\mathcal B_2(v,w)=-v\,\partial_x w
\tag{3}
\]

produces Burgers convection.

Time-independent nonhomogeneous boundary data are removed by writing

\[
u(x,t)=b(x)+v(x,t),
\tag{4}
\]

where \(b\) satisfies the boundary values; known lifting terms enter \(f\) and \(\mathcal F\).

All norm estimates below concern the finite-dimensional method-of-lines system, avoiding assumptions that unbounded PDE products map a Sobolev space into itself:

\[
\dot{\mathbf{u}}
=F_n(\mathbf{u})+\mathbf{f},
\qquad \mathbf{u}\in\mathbb C^n.
\tag{5}
\]

Ordinary symbols denote continuous scalar fields $u,W_j$, bold upright symbols semidiscrete vectors $\mathbf u,\mathbf W_j$, and bold italic symbols continuous vector fields $\boldsymbol U$. Uppercase letters denote matrices and calligraphic/script letters differential operators. Sections~2--3 use Euclidean norms and adjoints unless stated otherwise; local bounds may depend on $n$.

For systems and exact spatial reductions, $n$ denotes the active semidiscrete state dimension, not the number $N_x$ of points per spatial direction: use $n_{\rm raw}$ before reduction and $n_{\rm red}=n$ after it. Product-symmetry compression is a separate operation on the lifted state.

\subsubsection{Exact frozen--Fréchet decomposition}\label{exact-frozen-fruxe9chet-decomposition}

Let \(\mathcal U_n\subset\mathbb C^n\) be an admissible open neighborhood and let \(\bar{\mathbf{u}}\in\mathcal U_n\) be fixed in time, or fixed after a change to a co-moving coordinate. Define

\begin{align}
\mathbf{u}&=\bar{\mathbf{u}}+\mathbf{w},\qquad
\mathbf{w}(0)=\mathbf{w}_{\mathrm{in}}:=\mathbf{u}_{\mathrm{in}}-\bar{\mathbf{u}},
\tag{8}\\
\mathbf{r}_{\bar{\mathbf{u}}}&:=F_n(\bar{\mathbf{u}})+\mathbf{f},
\tag{9}\\
\mathcal R_{\bar{\mathbf{u}}}(\mathbf{w})
&:=F_n(\bar{\mathbf{u}}+\mathbf{w})-F_n(\bar{\mathbf{u}})
-DF_n(\bar{\mathbf{u}})\mathbf{w}.
\tag{10}
\end{align}

With the frozen Oseen matrix \(G_{\bar{\mathbf{u}}}:=DF_n(\bar{\mathbf{u}})\), Eq. (5) is exactly equivalent to

\[
\dot{\mathbf{w}}
=G_{\bar{\mathbf{u}}}\mathbf{w}+\mathbf{r}_{\bar{\mathbf{u}}}
+\mathcal R_{\bar{\mathbf{u}}}(\mathbf{w}),
\qquad \mathbf{w}(0)=\mathbf{w}_{\mathrm{in}}.
\tag{11}
\]

Equation (11) is an exact regrouping that retains the nonlinear remainder in the homotopy hierarchy.

Before discretization, the corresponding differential operator for Eq. (2) is

\[
\mathcal G_{\bar u}v
=\mathcal A v+
\sum_{r=2}^{d}\sum_{\ell=1}^{r}
\mathcal B_r(
\bar u,\ldots,\underset{\ell}{v},\ldots,\bar u).
\tag{12}
\]

\subsubsection{Second-order remainder bound}\label{second-order-remainder-bound}

\textbf{Theorem 1 (local quadratic bound for the semidiscrete frozen remainder).}\\
Let \(X_n=\mathbb C^n\), let \(\mathcal U_n\subset X_n\) be open, and let \(F_n\in C^2(\mathcal U_n;X_n)\). For this result, assume the local regularity bound

\[
\sup_{\mathbf{z}\in\mathcal U_n}\|D^2F_n(\mathbf{z})\|
\le C_{F,2}.
\tag{12a}
\]

Fix \(\bar{\mathbf{u}}\in\mathcal U_n\) and define \(\mathcal R_{\bar{\mathbf{u}}}(\mathbf{w}):=F_n(\bar{\mathbf{u}}+\mathbf{w})-F_n(\bar{\mathbf{u}})-DF_n(\bar{\mathbf{u}})\mathbf{w}\). If the segment \(\{\bar{\mathbf{u}}+\theta \mathbf{w}:0\le\theta\le1\}\) lies in \(\mathcal U_n\), then

\[
\|\mathcal R_{\bar{\mathbf{u}}}(\mathbf{w})\|_{X_n}
\le\frac{C_{F,2}}{2}\|\mathbf{w}\|_{X_n}^2.
\tag{13}
\]

If \(\bar{\mathbf{u}}+B_R(0)\subset\mathcal U_n\), then for \(\|\mathbf{w}\|_{X_n},\|\mathbf{v}\|_{X_n}\le R\),

\[
\|\mathcal R_{\bar{\mathbf{u}}}(\mathbf{w})-\mathcal R_{\bar{\mathbf{u}}}(\mathbf{v})\|_{X_n}
\le C_{F,2}R\|\mathbf{w}-\mathbf{v}\|_{X_n}.
\tag{14}
\]

\textbf{Proof.} Taylor's integral remainder gives Eq. (13); the mean-value formula gives Eq. (14). $\square$

If only the discretization \(A_n\) of the original linear part \(\mathcal A\) is used, define \(N_n(\mathbf{u}):=F_n(\mathbf{u})-A_n\mathbf{u}\). The correction equation then contains

\[
DN_n(\bar{\mathbf{u}})\mathbf{w}+\mathcal R_{\bar{\mathbf{u}}}(\mathbf{w}),
\tag{16}
\]

The first term is \(O(\|\mathbf{w}\|)\), whereas the frozen--Oseen remainder is \(O(\|\mathbf{w}\|^2)\). This local improvement translates to a homotopy advantage when the semigroup of \(G_{\bar{\mathbf{u}}}\) is sufficiently controlled.

\subsubsection{Mild solution and contraction criterion}\label{mild-solution-and-contraction-criterion}

For this analysis, require the frozen operator

\[
G_{\bar{\mathbf{u}}}:=DF_n(\bar{\mathbf{u}})
\tag{12b}
\]

to generate a strongly continuous semigroup or group \(S_{\bar{\mathbf{u}}}(t)\) satisfying

\[
\|S_{\bar{\mathbf{u}}}(t)\|
\le M_{\bar{\mathbf{u}}}e^{\omega_{\bar{\mathbf{u}}}t},
\qquad 0\le t\le T.
\tag{12c}
\]

The mild form of Eq. (11) is

\[
\mathbf{w}(t)=S_{\bar{\mathbf{u}}}(t)\mathbf{w}_{\mathrm{in}}
+\int_0^tS_{\bar{\mathbf{u}}}(t-s)
\left[\mathbf{r}_{\bar{\mathbf{u}}}+\mathcal R_{\bar{\mathbf{u}}}(\mathbf{w}(s))\right]\mathrm ds.
\tag{17}
\]

Define

\[
E_\omega(T)=e^{\max(\omega,0)T},
\qquad
\Phi_\omega(T)=
\begin{cases}
(e^{\omega T}-1)/\omega,&\omega\ne0,\\
T,&\omega=0.
\end{cases}.
\tag{18}
\]

\textbf{Theorem 2 (sufficient local existence and contraction condition).}\\
Assume the local second-derivative bound in Eq. (12a) and the frozen-semigroup estimate in Eq. (12c). Suppose there is \(R>0\) such that \(\bar{\mathbf{u}}+B_R(0)\subset\mathcal U_n\),

\[
\begin{aligned}
&M_{\bar{\mathbf{u}}}E_{\omega_{\bar{\mathbf{u}}}}(T)\|\mathbf{w}_{\mathrm{in}}\|_{X_n}\\
&\quad+
M_{\bar{\mathbf{u}}}\Phi_{\omega_{\bar{\mathbf{u}}}}(T)
\left(
\|\mathbf{r}_{\bar{\mathbf{u}}}\|_{X_n}
+\frac{C_{F,2}}{2}R^2
\right)
\le R
\end{aligned}
\tag{19}
\]

and

\[
q_{\bar{\mathbf{u}}}:=
M_{\bar{\mathbf{u}}}C_{F,2}R\Phi_{\omega_{\bar{\mathbf{u}}}}(T)<1.
\tag{20}
\]

Then the Duhamel map in Eq. (17) is a contraction on the closed radius-\(R\) ball of \(C([0,T];X_n)\), and Eq. (11) has a unique mild solution in that ball.

\textbf{Proof.} Equations (12c) and (13) imply the self-map bound (19). Equations (12c) and (14) imply that the Duhamel map has Lipschitz constant no larger than \(q_{\bar{\mathbf{u}}}\). The claim follows from the Banach fixed-point theorem. \(\square\)

\subsubsection{Conditional comparison with the original linear operator}\label{a-conditional-comparison-with-the-original-linear-operator}

Let \(M_A,\omega_A\) denote a semigroup bound for \(A_n\), and set \(L_1=\|DN_n(\bar{\mathbf{u}})\|\). Within a radius-\(R\) ball, the Lipschitz bounds for the two remainders are

\[
L_{\rm base}(R)\le L_1+C_{F,2}R,
\qquad
L_{\rm Oseen}(R)\le C_{F,2}R.
\tag{21}
\]

Consequently, a sufficient condition for a strictly smaller contraction upper bound is

\[
M_{\bar{\mathbf{u}}}\Phi_{\omega_{\bar{\mathbf{u}}}}(T)C_{F,2}R
<
M_A\Phi_{\omega_A}(T)(L_1+C_{F,2}R),
\tag{22}
\]

provided that the self-map condition holds for both decompositions.

Equation (22) gives a conditional comparison: the reduction in the nonlinear Lipschitz constant must outweigh any increase in transient amplification by the linear propagator.

\subsubsection{Homotopy deformation hierarchy}\label{homotopy-deformation-equations}

Let

\[
\mathscr L_{\bar{\mathbf{u}}}:=\partial_t-G_{\bar{\mathbf{u}}}.
\tag{23}
\]

Let \(\mathbf{W}_0\) denote the zeroth deformation coefficient. We use the standard homotopy deformation equation \cite{ref9,ref10,ref11}

\[
(1-q)\mathscr L_{\bar{\mathbf{u}}}[\boldsymbol{\Phi}(q)-\mathbf{W}_0]
=q\hbar H_{\rm HAM}\,
\mathscr N_{\bar{\mathbf{u}}}[\boldsymbol{\Phi}(q)],
\tag{24}
\]

where

\[
\mathscr N_{\bar{\mathbf{u}}}[\mathbf{w}]
=\mathscr L_{\bar{\mathbf{u}}}\mathbf{w}-\mathbf{r}_{\bar{\mathbf{u}}}
-\mathcal R_{\bar{\mathbf{u}}}(\mathbf{w}).
\tag{25}
\]

The embedding parameter is \(q\in[0,1]\), \(\hbar\ne0\) controls convergence, and \(H_{\rm HAM}\) is an auxiliary function. Expanding in the \(j\)-th deformation coefficients \(\mathbf{W}_j\),

\[
\boldsymbol{\Phi}(q)=\sum_{j=0}^{\infty}q^j\mathbf{W}_j
\tag{26}
\]

and taking the practical choice \(\hbar=-1\), \(H_{\rm HAM}=1\), with

\[
\mathscr L_{\bar{\mathbf{u}}}\mathbf{W}_0=\mathbf{r}_{\bar{\mathbf{u}}},
\qquad \mathbf{W}_0(0)=\mathbf{w}_{\mathrm{in}},
\tag{27}
\]

yields

\[
\boxed{
\mathscr L_{\bar{\mathbf{u}}}\mathbf{W}_j
=
[q^{j-1}]
\mathcal R_{\bar{\mathbf{u}}}
\left(\sum_{\ell=0}^{j-1}q^\ell \mathbf{W}_\ell\right),
\qquad \mathbf{W}_j(0)=0,\quad j\ge1.
}
\tag{28}
\]

For the polynomial remainder

\[
\mathcal R_{\bar{\mathbf{u}}}(\mathbf{w})
=\sum_{r=2}^{d}
\widetilde{\mathcal B}_{r,\bar{\mathbf{u}}}(\mathbf{w},\ldots,\mathbf{w}),
\tag{29}
\]

Eq. (28) becomes

\[
\mathscr L_{\bar{\mathbf{u}}}\mathbf{W}_j
=
\sum_{r=2}^{d}
\sum_{\substack{j_1+\cdots+j_r=j-1\\j_\ell\ge0}}
\widetilde{\mathcal B}_{r,\bar{\mathbf{u}}}
(\mathbf{W}_{j_1},\ldots,\mathbf{W}_{j_r}).
\tag{30}
\]

All orders use the same frozen operator. The order-\(m\) approximation is

\[
\mathbf{u}^{(m)}=\bar{\mathbf{u}}+\sum_{j=0}^{m}\mathbf{W}_j.
\tag{31}
\]

Let \(C_T:=C([0,T];X_n)\). If there exist an index \(j_0\ge0\) and a constant \(0<\alpha_{\rm HAM}<1\) such that

\[
\|\mathbf{W}_{j+1}\|_{C_T}
\le\alpha_{\rm HAM}\|\mathbf{W}_j\|_{C_T},
\qquad j\ge j_0,
\tag{32}
\]

then the deformation series converges absolutely in \(C_T\). In the semidiscrete polynomial setting, uniform convergence permits the limit to pass through the mild form of the finite-dimensional recursion. If \(\mathbf{w}_\infty:=\sum_{j=0}^{\infty}\mathbf{W}_j\) lies in the same closed radius-\(R\) ball used in Theorem 2, then it satisfies Eq. (17) at \(q=1\), and uniqueness in that ball identifies \(\bar{\mathbf{u}}+\mathbf{w}_\infty\) with \(\mathbf{u}\). Consequently, for every \(m\ge j_0\),

\[
\left\|\mathbf{u}-\mathbf{u}^{(m)}\right\|_{C_T}
\le\frac{\alpha_{\rm HAM}}{1-\alpha_{\rm HAM}}
\|\mathbf{W}_m\|_{C_T}.
\tag{33}
\]

Equation (33) requires an infinite-tail ratio. Mesh-uniform order additionally requires uniform control of that ratio and the first tail coefficient \(\|\mathbf W_{j_0}^{(n)}\|_{C_T}\).

For \(\hbar=-1\), \(H_{\rm HAM}=1\), and Eq. (27), the zero-order equation is equivalently

\[
\mathscr L_{\bar{\mathbf{u}}}\boldsymbol{\Phi}(q)
=\mathbf{r}_{\bar{\mathbf{u}}}+q\,\mathcal R_{\bar{\mathbf{u}}}(\boldsymbol{\Phi}(q)).
\tag{33a}
\]

\textbf{Proposition 1 (a sufficient condition for convergence at \(q=1\)).}\\
Assume that \(\mathcal R_{\bar{\mathbf{u}}}\) has a complex-analytic extension to the radius-\(R\) ball. If there is a \(\rho>1\) such that the self-map version of Eq. (19), with \(C_{F,2}R^2/2\) replaced by \(\rho C_{F,2}R^2/2\), holds and

\[
\rho M_{\bar{\mathbf{u}}}C_{F,2}R\Phi_{\omega_{\bar{\mathbf{u}}}}(T)<1,
\tag{33b}
\]

then the solution \(\boldsymbol{\Phi}(q)\) of Eq. (33a) is analytic for \(|q|\le\rho\). If \(\sup_{|q|=\rho}\|\boldsymbol{\Phi}(q)\|_{C_T}\le C_\rho\), Cauchy's estimate gives

\[
\left\|
\boldsymbol{\Phi}(1)-\sum_{j=0}^{m}\mathbf{W}_j
\right\|_{C_T}
\le
\frac{C_\rho\,\rho^{-(m+1)}}{1-\rho^{-1}}.
\tag{33c}
\]

\textbf{Proof.} For every \(|q|\le\rho\), Eqs. (13)--(14) and the stated inequalities make the Duhamel map corresponding to Eq. (33a) a uniform contraction on the complex radius-\(R\) ball. Analytic dependence of a uniformly contractive fixed point on \(q\) gives analyticity of \(\boldsymbol{\Phi}\). Applying Cauchy's coefficient estimate and summing the geometric tail yields Eq. (33c). \(\square\)

QHAM~\cite{ref12} retains the affine residual and closes non-iterative quadratic coefficient products. In its notation, set $N_0=\mathbf r$, $N_1(\mathbf w)=G\mathbf w$, $N_2=\mathcal R$, $L=\partial_t-G$, and $U_0=\mathbf W_0(t)$ satisfying $L\mathbf W_0=\mathbf r$. The prescribed physical profile is added at reconstruction. Our contribution is the auxiliary-dependent convergence and resource analysis built on that closure.

\paragraph{Auxiliary-dependent sufficient order}
For comparison with the numerical diagnostics, write $C_2:=C_{F,2}$ for the same Hessian bound.
Keep the profile, nonlinear semidiscrete system and interval fixed, and prescribe a finite set of candidate triples $(A,R,\rho)$ with constant auxiliary matrix $A$, $R>0$ and $\rho>1$. Write $\Delta_A=\|G-A\|$, $\|e^{tA}\|\le M_Ae^{\omega_A t}$, $b_A=M_A\Phi_{\omega_A}(T)$, and $a_A=M_AE_{\omega_A}(T)\|\mathbf w_{\rm in}\|+b_A\|\mathbf r\|$.
\textbf{Corollary 1.1 (sufficient auxiliary-order certificate).}
Assume that the remainder is analytic on a neighborhood of the complex radius-$R$ ball and satisfies Eqs. (13)--(14) there. If
\[
\begin{aligned}
 a_A+\rho b_A(\Delta_AR+C_{F,2}R^2/2)&\le R,\\
 \rho b_A(\Delta_A+C_{F,2}R)&<1,
\end{aligned}
\tag{33e}\label{eq:auxiliary-screen}
\]
then the embedding $(\partial_t-A)\boldsymbol\Phi_A=\mathbf r+q[(G-A)\boldsymbol\Phi_A+\mathcal R(\boldsymbol\Phi_A)]$, with initial value $\mathbf w_{\rm in}$, converges at $q=1$ to the same nonlinear solution, and
\[
\begin{aligned}
 \|\mathbf u-\mathbf u_A^{(m)}\|_{C_T}&\le\frac{R\rho^{-m}}{\rho-1},\\
 m_A^{\rm cert}&=\max\left\{0,\left\lceil
 \frac{\log[R/(\varepsilon_{\rm HAM}(\rho-1))]}{\log\rho}
 \right\rceil\right\}
\end{aligned}
\tag{33f}\label{eq:certified-order}
\]
is sufficient for an absolute semidiscrete homotopy tolerance $\varepsilon_{\rm HAM}>0$.
\textbf{Proof.} The modified remainder has ball bound $\Delta_AR+C_{F,2}R^2/2$ and Lipschitz bound $\Delta_A+C_{F,2}R$. Equation~\eqref{eq:auxiliary-screen} gives a uniform self-map and strict contraction for $|q|\le\rho$. Proposition 1's Cauchy argument applies with $C_\rho\le R$; solving the tail bound for $m$ gives Eq.~\eqref{eq:certified-order}. The added hierarchy term $(G-A)\mathbf W_{j-1}$ lowers the weight in Eq. (45) by $d-1$, so the same complete closure and dimension formulas remain valid. $\square$

Full linearization removes $\Delta_A$ but can increase $b_A$. For centered Burgers, frozen advection $A=\nu D_2-\operatorname{diag}(\bar{\mathbf u})D_1$ gives $G-A=-\operatorname{diag}(D_1\bar{\mathbf u})$ and $\Delta_A=\|D_1\bar{\mathbf u}\|_\infty$. A verified logarithmic-norm bound permits $M_A=1$ and $\omega_A\ge\mu_2(A)$; eigenvalues alone do not. An empty certified set is inconclusive, not evidence of divergence. At fixed $n,d$, the smallest certified order minimizes the complete ordered or symmetric dimension among the certified choices, not the actual required order or total quantum cost. Comparing the latter also requires the candidate-specific access, shift, normalization and output factors of Section~\ref{resource-analysis-and-conditional-scaling-comparison}.

\paragraph{Residual-based error tube}
For a differentiable approximation $\mathbf v$ with the correct initial value, define $\mathbf d=\dot{\mathbf v}-F_n(\mathbf v)-\mathbf f$. Suppose $\|\mathbf v-\bar{\mathbf u}\|\le S$, $\|\mathbf d\|\le d_*$, $\mu_2(G)\le\omega$, and the Hessian bound holds throughout the radius-$(S+E)$ neighborhood. Then
\[
 d_*\Phi_{\omega+C_{F,2}S+C_{F,2}E/2}(T)<E
 \quad\Longrightarrow\quad
 \|\mathbf u-\mathbf v\|_{C_T}<E.
\tag{33g}\label{eq:residual-tube}
\]
Indeed, $e=\|\mathbf u-\mathbf v\|$ obeys $D^+e\le(\omega+C_{F,2}S)e+C_{F,2}e^2/2+d_*$. Before a first crossing of $E$, linear comparison bounds $e$ by the left-hand side, excluding that crossing. For a computed trajectory, a certified $d_*$ must include its numerical evolution defect; residual samples alone are not upper bounds.

\subsubsection{Frozen-profile construction and selection}\label{selection-of-the-frozen-profile}

The tested choices are tied to the governing dynamics and fixed before LCHS propagation. Table~\ref{tab:tested-profiles} lists only profile constructions used in Section~\ref{numerical-experiments}; their numerical support is case-specific.

\begin{table}[H]
\centering
\caption{Frozen profiles used in the numerical experiments; evidence is limited to the cited tests.}
\label{tab:tested-profiles}
\TableFont
\begin{tabularx}{\linewidth}{@{}p{0.25\linewidth}X@{}}
\toprule
Profile & Construction and tested cases \\
\midrule
Initial state & $\bar{\mathbf u}=\mathbf u_{\rm in}$; unfitted Burgers, KdV, 2D Burgers/ZK, and vortex pair (Tables~\ref{tab:auxiliary-ablation}, \ref{tab:kdv-finite-k}, \ref{tab:extension-errors}; Fig.~\ref{fig:ns2d-vortex-pair}). \\
\addlinespace[2pt]
Layer enrichment & Multiplicative Burgers boundary-layer profile, Eq.~(90), calibrated in-sample (Table~\ref{tab:burgers-common-grid}). \\
\addlinespace[2pt]
Acoustic mean & Time average in Eq.~\eqref{eq:acoustic-profile} for 1D and 2D NS; a 1D endpoint average is also tested (Tables~\ref{tab:extension-errors}, \ref{tab:ns2d-errors}). \\
\bottomrule
\end{tabularx}
\end{table}

\Needspace{5\baselineskip}
Initial-state freezing removes the initial correction, while the Burgers layer enrichment and NS acoustic mean target their respective transport regimes. For the NS predictor variables $\mathbf a^a$ in Eqs.~(111) and~(114), excluding the derived reciprocal variable, the fixed mean is
\[
\bar{\mathbf a}_{\rm ac}=\frac1T\int_0^T\mathbf a^a(t)\,\mathrm dt.
\tag{34}\label{eq:acoustic-profile}
\]
KdV still freezes the initial wave rather than tracking it. These tests do not rank profiles across PDEs, and profile preparation enters end-to-end resource accounting.

\paragraph{Sensitivity to the selected profile}

The effect of profile accuracy on the frozen propagator can be quantified independently of a particular optimization algorithm.

\textbf{Proposition 2 (profile perturbation bound).}\\
Let \(\mathbf{p}\) and \(\mathbf{q}\) be two admissible semidiscrete profiles whose connecting segment lies in \(\mathcal U_n\), with \(\|D^2F_n\|\le C_{F,2}\) on that segment. Write \(G_{\mathbf{p}}=DF_n(\mathbf{p})\), \(G_{\mathbf{q}}=DF_n(\mathbf{q})\), and suppose \(\|e^{tG_{\mathbf{a}}}\|\le M_{\mathbf{a}} e^{\omega_{\mathbf{a}} t}\) for \(\mathbf{a}=\mathbf{p},\mathbf{q}\). Then

\[
\|G_{\mathbf{p}}-G_{\mathbf{q}}\|\le C_{F,2}\|\mathbf{p}-\mathbf{q}\|,
\qquad
|\mu_2(G_{\mathbf{p}})-\mu_2(G_{\mathbf{q}})|\le C_{F,2}\|\mathbf{p}-\mathbf{q}\|,
\tag{35a}
\]

and

\[
\|e^{tG_{\mathbf{p}}}-e^{tG_{\mathbf{q}}}\|
\le M_{\mathbf{p}}M_{\mathbf{q}} C_{F,2}\|\mathbf{p}-\mathbf{q}\|
\int_0^t e^{\omega_{\mathbf{p}}(t-s)+\omega_{\mathbf{q}}s}\,\mathrm ds.
\tag{35b}
\]

With a common bound \(M_{\mathbf{p}},M_{\mathbf{q}}\le M\), \(\omega_{\mathbf{p}},\omega_{\mathbf{q}}\le\omega\), the right-hand side is at most \(M^2 C_{F,2}\|\mathbf{p}-\mathbf{q}\|\,t e^{\omega t}\).

\textbf{Proof.} Integrating \(D^2F_n\) along the segment gives the first bound. The largest eigenvalue of a Hermitian matrix is Lipschitz in the operator norm; applying this fact to the Hermitian parts of \(G_{\mathbf{p}}\) and \(G_{\mathbf{q}}\) gives the logarithmic-norm bound. Finally,

\[
e^{tG_{\mathbf{p}}}-e^{tG_{\mathbf{q}}}
=\int_0^t e^{(t-s)G_{\mathbf{p}}}(G_{\mathbf{p}}-G_{\mathbf{q}})e^{sG_{\mathbf{q}}}\,\mathrm ds,
\tag{35c}
\]

and taking norms proves Eq. (35b). \(\square\)

Consequently, fitting a profile accurately in a norm that controls \(C_{F,2}\|\mathbf{p}-\mathbf{q}\|\) can preserve an already favorable frozen propagator, but a small field misfit alone is insufficient when the semidiscrete derivative constants grow under refinement. For a candidate \(\mathbf{p}\), Eq. (35a) also supplies a robustness margin: if \(\mu_2(G_{\mathbf{p}})\le-\gamma_0<0\), then every profile \(\mathbf{q}\) satisfying \(C_{F,2}\|\mathbf{p}-\mathbf{q}\|<\gamma_0\) retains a negative logarithmic norm.

For fixed \(F_n\), consistent recomputation of the residual and remainder preserves the exact solution. Changing the profile affects the finite-order tail, not the physical model. Profile-representation and assembly errors must therefore be distinguished from that tail; report residual and stability metrics alongside field-fit errors.

\subsection{Quantum-compatible linearization by product lifting}\label{quantum-compatible-linearization-by-product-lifting}

Following the coefficient-product lifting of Ref. \cite{ref12}, we close the prescribed finite homotopy hierarchy. The remaining approximation is the homotopy tail, not an additional truncation of the lifted system.

\subsubsection{Finite closure of the quadratic hierarchy}\label{quadratic-hierarchy}

The spatially discretized quadratic hierarchy is

\[
\dot{\mathbf{W}}_j
=G\mathbf{W}_j+\delta_{j0}\mathbf{r}
+\mathbf 1_{\{j\ge1\}}
\sum_{a+b=j-1}B(\mathbf{W}_a,\mathbf{W}_b),
\qquad 0\le j\le m,
\tag{36}
\]

Here \(G\in\mathbb C^{n\times n}\), \(\mathbf{r}\in\mathbb C^n\), and \(B\) are the fixed Oseen matrix, residual, and bilinear map; \(\delta_{j0}\) and \(\mathbf 1_{\{j\ge1\}}\) are the Kronecker delta and indicator.

For an ordered multi-index \(\boldsymbol\alpha\in\mathbb N_0^p\), with \(|\boldsymbol\alpha|=\sum_\ell\alpha_\ell\), retain

\[
\mathbf{Y}_{\boldsymbol\alpha}
=\mathbf{W}_{\alpha_1}\otimes\cdots\otimes \mathbf{W}_{\alpha_p},
\qquad
|\boldsymbol\alpha|+p\le m+1.
\tag{37}
\]

Differentiation preserves the weight $|\boldsymbol\alpha|+p$: replacing $\mathbf W_j$ by $\mathbf W_a\otimes\mathbf W_b$, $a+b=j-1$, leaves it unchanged, while residual insertion lowers it. Hence Eq. (37) closes exactly.

\subsubsection{Lifted dimension and target recovery}\label{exact-block-and-dimension-count-for-the-quadratic-construction}

Append to the product state \(\mathbf{Y}\) the target block

\[
\mathbf{Y}_{-1}:=\sum_{j=0}^{m}\mathbf{W}_j.
\tag{42}
\]

The coordinate projector \(P_{\rm tar}\) selects \(\mathbf{Y}_{-1}\), and the physical approximation is \(\bar{\mathbf{u}}+P_{\rm tar}\mathbf{Y}\).

The resulting ordered, uniformly padded dimension is \cite{ref12}

\[
\begin{aligned}
D_m
&=n+\sum_{p=1}^{m+1}\binom{m+1}{p}n^p\\
&=(n+1)^{m+1}+n-1.
\end{aligned}
\tag{44}
\]

An exact symmetry reduction may compress this representation if initialization, dynamics, and target recovery are preserved. Discarding required product variables is instead an additional approximation.

\subsubsection{Extension to general polynomial nonlinearities}\label{polynomial-nonlinearities-of-general-degree}

For the degree-\(d\) hierarchy, assign

\[
\operatorname{wt}_d(\boldsymbol\alpha)
=(d-1)|\boldsymbol\alpha|+p.
\tag{45}
\]

Replacing \(\mathbf{W}_j\) by \(r\) factors with \(2\le r\le d\) and \(\sum_{\ell=1}^rj_\ell=j-1\) gives

\[
(d-1)\sum_{\ell=1}^{r}j_\ell+r
=(d-1)j+r-d+1
\le(d-1)j+1.
\tag{46}
\]

\textbf{Theorem 3 (finite polynomial closure).}\\
Fix a homotopy order \(m\ge0\) and a polynomial degree \(d\ge2\). Define the index set

\[
\mathcal I_{m,d}
=\left\{
\boldsymbol\alpha\in\mathbb N_0^p:
1\le p\le1+(d-1)m,\ 
\operatorname{wt}_d(\boldsymbol\alpha)
\le(d-1)m+1
\right\}
\tag{47}
\]

The variables indexed by \(\mathcal I_{m,d}\) are closed under the product-rule dynamics of Eq. (30), with an affine source or one constant coordinate. Their maximum tensor degree satisfies

\[
p_{\max}\le1+(d-1)m.
\tag{48}
\]

\textbf{Proof.} Linear terms preserve weight; Eq. (46) shows that nonlinear replacements cannot increase it. Replacing \(\mathbf{W}_0\) by its affine source lowers the degree. Thus every product derivative is an affine linear combination within \(\mathcal I_{m,d}\). \(\square\)

The retained products therefore satisfy

\[
\boxed{
\dot{\mathbf{Y}}(t)=A_m\mathbf{Y}(t)+\mathbf{b}_m,
\qquad \mathbf{Y}(0)=\mathbf{Y}_{\rm in}.
}
\tag{49}
\]

where \(A_m\) and \(\mathbf{b}_m\) remain fixed on \([0,T]\); propagation requires no nonlinear products.

\textbf{Corollary 3.1 (equivalence with the truncated deformation hierarchy).}\\
Initialize every lifted block consistently:

\[
\mathbf{Y}_{\boldsymbol\alpha}(0)
=\bigotimes_{\ell=1}^{p}\mathbf{W}_{\alpha_\ell}(0),
\qquad
\mathbf{Y}_{-1}(0)=\mathbf{w}_{\rm in}.
\tag{49a}
\]

Only the all-zero-index products \(\mathbf{w}_{\rm in}^{\otimes p}\) are initially nonzero. Exact propagation preserves

\[
\mathbf{Y}_{\boldsymbol\alpha}(t)
=\bigotimes_{\ell=1}^{p}\mathbf{W}_{\alpha_\ell}(t),
\qquad
\mathbf{Y}_{-1}(t)=\sum_{j=0}^{m}\mathbf{W}_j(t).
\tag{49b}
\]

\textbf{Proof.} Products formed from the unique triangular deformation solution satisfy the same affine system and initial data. Uniqueness proves the first identity; summing the deformation equations proves the target identity. \(\square\)

Thus exact lifting introduces no truncation beyond order \(m\). Violations of Eq. (49b) diagnose initialization or propagation errors, separately from the HAM tail.

For general degree \(d\), define

\[
J_p=
\left\lfloor
\frac{(d-1)m+1-p}{d-1}
\right\rfloor ,
\qquad
1\le p\le1+(d-1)m.
\tag{50a}
\]

The degree-\(p\) block count and total dimensions are

\[
N_p=\binom{J_p+p}{p}.
\tag{50b}
\]

\[
D_{m,d}^{\rm lift}
=\sum_{p=1}^{1+(d-1)m}N_pn^p,
\qquad
D_{m,d}=D_{m,d}^{\rm lift}+n
\tag{50c}
\]

including the target block. For \(d=2\), this recovers Eq. (44). Below, \(D_m:=D_{m,d}\) for fixed \(d\); its polynomial degree in \(n\) grows with \(m\).

\textbf{Lemma 1 (exact product-symmetry quotient).}\\
Set \(x_{j,i}=(\mathbf{W}_j)_i\) and give each tagged scalar variable \((j,i)\) weight \(1+(d-1)j\). Let \(\mathcal M\) contain every nonempty multiset \(\mu\) of these variables with total weight at most \(B=1+(d-1)m\). If \(\nu_a(\mu)\) is the multiplicity of tag \(a\), set \(h_\mu=|\mu|!/\prod_a\nu_a(\mu)!\) and \(Y_\mu^{\rm sym}=\sqrt{h_\mu}\prod_a x_a^{\nu_a(\mu)}\). Let the isometry \(S\) put \(Y_\mu^{\rm sym}/\sqrt{h_\mu}\) in each of the \(h_\mu\) ordered coordinates and act identically on the target block. For the complete ordered product-rule generator in Eq. (49),
\[
S^*S=I,\qquad A_mS=SB_m^{\rm sym},\qquad B_m^{\rm sym}=S^*A_mS.
\tag{50c.1}\label{eq:symmetric-intertwining}
\]
If \(\mathbf{b}_m=S\boldsymbol{\beta}_m\) and \(\mathbf{Y}_{\rm in}=S\mathbf{Y}^{\rm sym}_{\rm in}\), the exact trajectories satisfy \(\mathbf{Y}(t)=S\mathbf{Y}^{\rm sym}(t)\), and their target fields coincide. In particular,
\[
D_{m,d}^{\rm sym}=n+\sum_{\ell=1}^{B}[z^\ell]
\prod_{j=0}^{m}(1-z^{1+(d-1)j})^{-n}.
\tag{50c.2}\label{eq:symmetric-dimension}
\]
\textbf{Proof.} Write every scalar coefficient equation as \(\dot x_a=\sum_\gamma c_{a,\gamma}x^\gamma\), adding contributions from ordered, possibly nonsymmetric multilinear maps that yield the same commutative monomial. The product rule gives
\[
\dot Y_\mu^{\rm sym}=\sum_{a,\gamma}\nu_a(\mu)c_{a,\gamma}
\sqrt{\frac{h_\mu}{h_{\mu-e_a+\gamma}}}
Y_{\mu-e_a+\gamma}^{\rm sym},\qquad Y_\varnothing^{\rm sym}=h_\varnothing=1.
\tag{50c.3}
\]
Equation (46) keeps every nonconstant term within \(\mathcal M\). The right-hand side depends only on the multiset, so the ordered generator preserves \({\rm im}\,S\) and satisfies Eq. (50c.1). Compatible initialization and uniqueness prove the trajectory and target identities. \(\square\)

Here $\mathbf Y^{\rm sym}$ denotes unscaled quotient coordinates, distinct from the encoded state $\mathbf Z=D^{-1}\mathbf Y$ introduced below. This equivalence concerns exact affine dynamics. A finite LCHS rule on \(B_m^{\rm sym}\) need not equal the projection of the same rule on \(A_m\), because \(A_m^*\) need not preserve \({\rm im}\,S\). The experiments apply the finite rule to a similarity-scaled compressed generator, map its output back, and compare it with exact compressed propagation. The isometry preserves this unscaled-coordinate error after embedding both states in the ordered space; it does not preserve norms in the similarity-scaled coordinates.

The integer \(D_m\) need not itself be a power of two. The following padding applies equally to the ordered system or its exact symmetric quotient, with the corresponding active dimension, generator, source, initialization, and target map. For an amplitude-encoded implementation, set \(q_m=\lceil\log_2D_m\rceil\), \(\widehat D_m=2^{q_m}\), and let \(J:\mathbb C^{D_m}\to\mathbb C^{\widehat D_m}\) append zero coordinates. We use the exact block embedding
\[
\widehat A_m=A_m\oplus0_{\widehat D_m-D_m},\qquad
\widehat{\mathbf{b}}_m=J\mathbf{b}_m,\qquad
\widehat{\mathbf{Y}}_{\rm in}=J\mathbf{Y}_{\rm in},\qquad
\widehat P_{\rm tar}=P_{\rm tar}J^*.
\tag{50d}\label{eq:padding}
\]
The inactive coordinates remain zero, so \(J^*\widehat{\mathbf{Y}}(t)=\mathbf{Y}(t)\); the homogeneous and source-channel finite-LCHS rules likewise preserve this subspace because every simulated Hermitian matrix is block diagonal, provided the spectral shift and quadrature parameters are retained from the active block. If the source is implemented by one affine coordinate, it is appended before padding and requires \(\lceil\log_2(D_m+1)\rceil\) system qubits. Here $q_{\rm state}:=q_m$ counts only state-register qubits, and $N_Q:=\widehat D_m$ is the padded dimension; for compression, replace $D_m$ by the actually encoded active dimension $D_{\rm act}$. The padding factor is below two, and no spatial grid is coarsened solely to make $D_m$ a power of two. The reported classical experiments act on the unpadded active block; Eq. (50d) specifies the mathematically equivalent register embedding, not an implemented quantum circuit.

\subsubsection{Sparse-access requirements}\label{sparsity}

Suppose the spatial matrix has row sparsity \(\sigma\), and each multilinear contraction has computable \(O(1)\) or \(\operatorname{polylog}n\) row access. Degree-\(p\) Kronecker sums have sparsity \(O(p\sigma)\). With \(s\) local bilinear channels, the quadratic construction gives \cite{ref12}

\[
s_m=O(m\sigma+s\,m^2).
\tag{51}
\]

where \(s_m\) is the lifted row sparsity. Implicit access avoids explicit enumeration of \(D_m\); dense spectral or nonlocal operators instead require a structured block encoding.

LCHS also uses \(A_m^\dagger\), so a sparse-Hamiltonian implementation needs efficient row and column access. A dense residual in \(I\otimes \mathbf{r}\) can destroy column sparsity; structured source encoding must then include its preparation cost and normalization \(\|\mathbf{r}\|\). Equation (51) alone does not bound these costs.

A centered-difference, constant-kinematic-viscosity NS model admits explicit ordered row/column access and finite input preparation. Dense residual insertion creates $O(Pn)$ reverse candidates at maximum degree $P=m+1$; a structured alternative charges $P\|\mathbf r\|/\eta$. This local ordered construction does not certify the Fourier, symmetry-compressed production layouts.

\subsection{LCHS realization of the lifted linear dynamics}\label{lchs-realization-of-the-lifted-linear-dynamics}

\subsubsection{Homogeneous and source propagation}\label{homogeneous-and-source-propagation}\label{discrete-operators}

Set \(A=-A_m\) and \(\mathbf{b}=\mathbf{b}_m\). Then

\[
\dot{\mathbf{Y}}=-A\mathbf{Y}+\mathbf{b}.
\tag{52}
\]

with solution

\[
\mathbf{Y}(t)=e^{-At}\mathbf{Y}_{\rm in}
+\int_0^t e^{-A(t-s)}\mathbf{b}\,\mathrm ds.
\tag{53}
\]

Direct exponential-action methods \cite{ref50,ref51,ref52,ref53,ref54,ref55,ref56} provide verification diagnostics, not LCHS experiment curves.

Let

\[
A=L+iH,\qquad
L=\frac{A+A^\dagger}{2},\qquad
H=\frac{A-A^\dagger}{2i}.
\tag{54}
\]

When \(L\succeq0\), the exact Cauchy-kernel identity is \cite{ref15}

\[
e^{-At}\mathbf{v}
=\int_{\mathbb R}g_0(k)\,
e^{-i(H+kL)t}\mathbf{v}\,\mathrm dk,
\qquad
g_0(k)=\frac{1}{\pi(1+k^2)}.
\tag{55}
\]

Generalized kernels require the hypotheses of Ref. \cite{ref16}. A finite rule is

\[
\mathcal Q_{K,M}(t)\mathbf{v}
=\sum_{\ell=1}^{M}
\omega_\ell^{(K)}
e^{-i(H+k_\ell L)t}\mathbf{v},
\qquad |k_\ell|\le K.
\tag{56}
\]

where \(K\) is the cutoff half-width and \(M\) the node count.

The experiments use the Gaussian-damped Cauchy kernel

\[
\begin{aligned}
g_{\gamma_{\rm ker},c_{\rm ker}}(k)
&=
\frac{\exp\!\left(
c_{\rm ker}-(k^2+1)/(4\gamma_{\rm ker}^2)-ic_{\rm ker}k
\right)}{\pi(1+k^2)},\\
\gamma_{\rm ker}
&=c_{\rm ker}^{-1}
\sqrt{c_{\rm ker}+\log[(1+1/(2\pi))/\varepsilon_{\rm ker}]}.
\end{aligned}
\tag{57}
\]

We discretize the \(k\)-integral by the equally spaced finite sum of Low and Somma~\cite{ref69}. For \(M=2J+1\), \(h=K/J\), and \(k_j=jh\), its coefficients and propagator are

\[
\omega_j=h\,g_{\gamma_{\rm ker},c_{\rm ker}}(k_j),\qquad
\mathcal Q_{K,M}(t)=\sum_{j=-J}^{J}\omega_j e^{-i(H+k_jL)t}.
\tag{57a}
\]

Both endpoints carry the full weight \(h\), and the complex coefficient sum is not divided out. This is the finite truncation of the uniform trapezoidal sum in Ref.~\cite{ref69}, not the half-endpoint trapezoid rule on a bounded interval. For \(K=2c_{\rm ker}\gamma_{\rm ker}^2\) and sufficiently small \(h\), that reference bounds the raw finite-sum operator error for \(L\succeq0\), without requiring \(H\) and \(L\) to commute. Numerical parameters are stated in Section~\ref{numerical-experiments}.

\subsubsection{\texorpdfstring{Separate \(K_1\) and \(K_2\)}{Separate K\_1 and K\_2}}\label{separate-k_1-and-k_2}

Separate cutoffs and node counts control homogeneous and source propagation:

\[
e^{-At}\mathbf{Y}_{\rm in}
\approx\mathcal Q_{K_1,M_1}(t)\mathbf{Y}_{\rm in},
\tag{58}
\]

\[
\int_0^te^{-A(t-s)}\mathbf{b}\,\mathrm ds
\approx
\mathcal S_{K_2,M_2}(t;\mathbf{b}).
\tag{59}
\]

The emulator integrates each source node by one-coordinate affine augmentation. Quantum access to that source map is assumed in Section~\ref{resource-analysis-and-conditional-scaling-comparison}. Over \(N_t\) intervals, \(K_1\) stays fixed while \(K_2(t_j)\) may vary.

Source accuracy requires a uniform propagator bound over the integration interval (Eq. (59b)); sufficient Low--Somma choices are stated in Eq. (68a). Joint cutoff/node refinement is essential because interval refinement alone can approach the wrong generator (Proposition C.2).

\subsubsection{Similarity scaling and spectral shifting}\label{similarity-scaling-and-spectral-shifting}

To balance tensor blocks, use

\[
\mathbf{Y}=D\mathbf{Z},\qquad
D|_{\text{degree }p}=\eta^{p-1}I,
\qquad D|_{\mathbf{Y}_{-1}}=D|_{p=1}=I
\tag{62}
\]

which maps \(A\) and \(\mathbf{b}\) to \(D^{-1}AD\) and \(D^{-1}\mathbf{b}\), preserving the exact physical solution and eigenvalues.

For \(\rho_0=\|\mathbf{w}_{\rm in}\|>0\), \(P=1+(d-1)m\), and consistent initialization,

\[
\|\mathbf{Z}_{\rm in}\|^2
=\rho_0^2\left[
1+\sum_{p=1}^{P}(\rho_0/\eta)^{2(p-1)}
\right],
\qquad
p_{\rm tar,in}
=\left[
1+\sum_{p=1}^{P}(\rho_0/\eta)^{2(p-1)}
\right]^{-1}.
\tag{62a}
\]

Thus \(\rho_0/\eta>1\) yields exponentially small initial target weight as \(P\) grows; balancing the matrix can worsen state preparation. This initial probability does not bound the final one.

If needed, make the Hermitian part positive semidefinite using

\[
L_{\delta_{\rm sh}}=L+\delta_{\rm sh} I,\qquad
\delta_{\rm sh}\ge\max\{0,-\lambda_{\min}(L)\}
\tag{63}
\]

and compensate by

\[
e^{-At}=e^{\delta_{\rm sh}t}e^{-(A+\delta_{\rm sh} I)t},
\tag{63a}
\]

The amplification \(e^{\delta_{\rm sh}t}\) enters normalization, conditioning and success probability. The emulator chooses its shift from a numerical estimate of \(\lambda_{\min}(L)\) and checks finite sums against direct lifted propagation. An operator-level LCHS guarantee additionally requires a certified lower spectral bound, with its shift cost charged in Eq. (68).

\subsubsection{Error budget}\label{error-budget}

Finite-rule operator defects are weighted by the interval-start state, source norm and subsequent propagation. The experiments report physical-field, spatial and complete-lift errors separately.

\Needspace{7\baselineskip}
\subsubsection{Non-iterative workflow}\label{non-iterative-workflow}

\begin{algorithm}[!htb]
\caption{Non-iterative Frozen--Oseen HAM-LCHS}
\label{alg:frozen-oseen-lchs}
\DontPrintSemicolon
\KwIn{\(F_n,\mathbf{f},\mathbf{u}_{\rm in}\), order \(m\), profile family, time partition \(\{t_j\}_{j=0}^{N_t}\), and tolerances}
\KwOut{Requested observables and verified error/resource diagnostics}
Homogenize boundaries and choose a declared fixed \(\bar{\mathbf{u}}\)\;
Freeze \(G=DF_n(\bar{\mathbf{u}})\), \(\mathbf{r}=F_n(\bar{\mathbf{u}})+\mathbf{f}\), and \(\mathbf{w}_{\rm in}=\mathbf{u}_{\rm in}-\bar{\mathbf{u}}\)\;
Construct consistent products and verify access to \(A_m,\mathbf{b}_m,\mathbf{Y}_{\rm in},P_{\rm tar}\) in the actually encoded layout under the stated access model\;
Choose \(D\), the spectral shift, \(K_1,M_1\), and \(\{K_2(t_j),M_2(t_j)\}\)\;
Build the scaled, shift-compensated rules; set \(\mathbf{Z}_0=D^{-1}\mathbf{Y}_{\rm in}\), \(\widetilde{\mathbf{b}}=D^{-1}\mathbf{b}_m\)\;
\For{\(j=0,\ldots,N_t-1\)}{
  \(\Delta t_j\leftarrow t_{j+1}-t_j\)\;
  \(\mathbf{Z}_{j+1}\leftarrow
  \mathcal Q_{K_1,M_1}(\Delta t_j)\mathbf{Z}_j+
  \mathcal S_{K_2(t_j),M_2(t_j)}(\Delta t_j;\widetilde{\mathbf{b}})\)\;
}
Recover requested observables from \(\mathbf{u}^{(m)}(T)=\bar{\mathbf{u}}+P_{\rm tar}D\mathbf{Z}_{N_t}\)\;
Check product identities and independent spatial, HAM, and finite-rule refinements\;
Report parameters, error contributions in Eq. (64), and resource costs\;
\end{algorithm}

The profile and lifted matrix remain fixed throughout; the time loop composes linear propagation, not iterative HAM updates.

\section{Resource analysis and conditional scaling comparison}\label{resource-analysis-and-conditional-scaling-comparison}

\subsection{Explicit classical lifted evolution}\label{explicit-classical-lifted-evolution}

For fixed \(m,d\), Eq. (50c) gives \(D_m=\Theta(n^{p_{\max}})\), \(p_{\max}=1+(d-1)m\). With row sparsity \(s_m\), explicit matrix storage, materialized-state storage and \(N_t\)-step work with \(r_{\rm mv}\) matrix--vector products per step satisfy

\begin{align}
S_C^{\rm matrix}&=O(s_mD_m), \tag{65}\\
S_C^{\rm state}&=\Omega(D_m), \tag{66}\\
T_C^{\rm lift}&=\widetilde O(N_t r_{\rm mv}s_mD_m). \tag{67}
\end{align}

Profile construction, matrix assembly and reference-solution generation are additional costs beyond Eq. (67). Writing the full lifted state takes \(\Omega(D_m)\) time, whereas few-observable tasks have different classical baselines. Where exact symmetry reduction is used, both classical and quantum lifted costs use \(D_{\rm sym}\) and its actual sparsity.

The same-output classical comparator may solve the original nonlinear PDE directly: local discretizations can take \(O(n)\) work per explicit step, and FFT-based methods can take \(O(n\log n)\), subject to stability and accuracy requirements. Low-rank or symmetric tensor implementations may also avoid materializing all \(D_m\) coordinates. Equation (67) measures explicit-lift work.

\subsection{Quantum access model}\label{quantum-access-model}

All quantum normalizations below use the encoded coordinates $\mathbf Z=D^{-1}\mathbf Y$, with $A_D=-D^{-1}A_mD=L+iH$, $\mathbf Z_0=D^{-1}\mathbf Y_{\rm in}$ and $\mathbf b_D=D^{-1}\mathbf b_m$. The degree scaling in Eq. (62) preserves the target; general decoding is stated separately. A block encoding of a matrix $B$ is a unitary $U_B$ whose success subblock satisfies
\[
 (\langle0^a|\otimes I)U_B(|0^a\rangle\otimes I)=B/\alpha_B,
 \qquad \alpha_B\ge\|B\|.
\tag{67b}\label{eq:block-encoding-definition}
\]
Matrix normalization $\alpha_B$, the outer LCU normalization $\mathcal B$, and the target signal norm are distinct quantities.

Quantum PDE costs depend on input access and the requested output as well as dimension \cite{ref30,ref31}. The following circuit accounting concerns a \emph{finite unitary-sum} realization; it is not a gate count for the classical affine-augmentation exponential used by the emulator in Eq. (59). Let \(\alpha_H,\alpha_L\) be block-encoding normalizations of \(H\) and \(L_{\delta_{\rm sh}}=L+\delta_{\rm sh}I\succeq0\), respectively. Let \(C_{\rm be}\) be the cost of one such oracle call, \(C_{\rm prep}\) the cost of preparing normalized \(\mathbf{Z}_0\) and \(\mathbf{b}_D\), and \(C_{\rm recon}\) the additional physical-field preparation cost. The matrix, source, coefficient, and target-projector oracles must all act without enumerating \(D_m\) entries. Write \(P_{\rm tar}\) for the orthogonal projector onto the block in Eq. (42), and \(\varepsilon_{\rm alg}\) for the coherent error budget before measurement.

The target block is the correction \(\mathbf{w}^{(m)}=P_{\rm tar}\mathbf{Z}\), whereas the physical field is \(\mathbf{u}^{(m)}=\bar{\mathbf{u}}+\mathbf{w}^{(m)}\). The normalized state \(|\mathbf{w}^{(m)}\rangle\) alone does not determine its amplitude or the sum with \(\bar{\mathbf{u}}\). Staged reconstruction requires the correction norm, phase-consistent profile preparation and linear combination of unitaries (LCU) addition, with normalization-to-signal ratio

\[
\Gamma_{\rm recon}
=\frac{\|\bar{\mathbf{u}}\|+\|\mathbf{w}^{(m)}\|}
{\|\bar{\mathbf{u}}+\mathbf{w}^{(m)}\|}.
\tag{67a}
\]

Equation (67a) applies to staged reconstruction from an already prepared correction state. Alternatively, the known profile can enter the original LCU as a third branch, avoiding intermediate correction-state postselection. Dense classical reconstruction remains an \(n\)-component output task.

The time-independent lifted system permits direct evaluation of Eq. (53) at a requested time \(t\); interval composition is needed for the reported emulator but is not intrinsic to the affine equation. To specify a quantum circuit, replace the source-node affine exponential by a finite Duhamel time sum of unitaries. Let the homogeneous and source \(k\)-rules have coefficient one-norms \(\lambda_i=\sum_{\ell=1}^{M_i}|\omega_\ell^{(i)}|\), \(i=1,2\), and let \(R_2\) be the number of source-time nodes. For positive time-bin weights summing to \(\Phi_{\delta_{\rm sh}}(t)\), an outer LCU has the explicit normalization and target success probability

\[
\begin{aligned}
\mathcal B(t)&=
e^{\delta_{\rm sh}t}\lambda_1\|\mathbf{Z}_0\|
+\Phi_{\delta_{\rm sh}}(t)\lambda_2\|\mathbf{b}_D\|,\\
p_{\rm sig}(t)&=
\frac{\|P_{\rm tar}\widetilde{\mathbf{Z}}(t)\|^2}{\mathcal B(t)^2}.
\end{aligned}
\tag{68}
\]

Equation (68) uses an orthogonal target projector in the encoded coordinates, with inputs transformed consistently. A general decoding map \(F=P_{\rm tar}D\) requires a block encoding \(F/\beta_F\), where \(\beta_F\ge\|F\|\); its success probability is \(\|F\widetilde{\mathbf{Z}}\|^2/(\beta_F\mathcal B)^2\).

Low and Somma~\cite{ref69} give an \emph{a priori} sufficient prescription for this rule. Restrict to \(c_{\rm ker}>0\), \(0<\varepsilon_{{\rm ker},i}\le0.9\), and \(0<\varepsilon_{{\rm q},i}\le4/15\), within their Theorems 2--3 ranges. For \(L\succeq0\), choose \(\gamma_i\) by Eq. (57) and set

\[
\begin{aligned}
K_i&=2c_{\rm ker}\gamma_i^2,\qquad M_i=2K_i/h_i+1,\\
h_i&\le
\frac{\pi}{
\frac12 t\|L\|+
\log\!\left[64e^{3c_{\rm ker}/2}/(15\varepsilon_{{\rm q},i})\right]},
\qquad K_i/h_i\in\mathbb N.
\end{aligned}
\tag{68a}
\]

For admissible error parameters, Eq. (68a) bounds the homogeneous operator defect by \(\varepsilon_{{\rm ker},1}+\varepsilon_{{\rm q},1}\). By Eq. (59b), a uniform integrand bound \(\varepsilon_{{\rm ker},2}+\varepsilon_{{\rm q},2}\) on \(0\le\tau\le t\) bounds the source-map defect by \(t(\varepsilon_{{\rm ker},2}+\varepsilon_{{\rm q},2})\), before source-time discretization. At fixed \(c_{\rm ker}\) and \(t\|L\|\), the sufficient choices have \(K_i=O(\log(1/\eta_i))\) and \(M_i=O([t\|L\|+\log(1/\eta_i)]\log(1/\eta_i))\), where \(\eta_i\) is the allocated integrand error. A shift replaces \(L\) by \(L_{\delta_{\rm sh}}\) and multiplies the corresponding defects by up to \(e^{\delta_{\rm sh}t}\). The experimental \(K=32\), \(M=385\) are tested through full-lift errors; they are not asserted to meet Eq. (68a) at \(\varepsilon_{\rm ker}=10^{-8}\).

For target-preserving scaling, adding the known profile directly to the LCU gives
\[
\mathcal B_{\rm phys}=\mathcal B+\|\bar{\mathbf{u}}\|,
\qquad
p_{\rm phys}=\frac{\|\bar{\mathbf{u}}+P_{\rm tar}\widetilde{\mathbf{Z}}\|^2}
{\mathcal B_{\rm phys}^2}.
\tag{68b}\label{eq:direct-profile-lcu}
\]
This requires phase-consistent profile preparation, but no intermediate estimate of the correction norm.

For the unitary source-time sum, a midpoint bin width \(h_t\) adds at most \(\tfrac12h_t\lambda_2\Omega_{K_2}\Phi_{\delta_{\rm sh}}(t)\) to the source-map defect, where \(\Omega_{K_2}=\|H\|+K_2\|L_{\delta_{\rm sh}}\|\). Hence an error allocation \(\varepsilon_t\) is ensured by \(R_2=O(1+t\lambda_2\Omega_{K_2}\Phi_{\delta_{\rm sh}}(t)/\varepsilon_t)\) equal-width bins. This is a sufficient construction, not the cost of the affine-augmentation emulator.

Let \(C_{\rm coeff}(M_1,M_2,R_2)\) include coherent preparation of the complex coefficient magnitudes and phases and \(C_{\rm sel}\) include their controlled selection. A block-encoded Hamiltonian simulation of any selected node has query cost \(\widetilde O((\alpha_H+K_{\max}\alpha_L)t+\log(1/\varepsilon_{\rm sim}))\), where \(K_{\max}=\max(K_1,K_2)\), under the usual controlled-simulation access assumptions \cite{ref69}. With a known positive lower bound on \(p_{\rm sig}\), coherent amplitude amplification gives the following \emph{one-output-time, conditional} gate-accounting bound:

\[
\boxed{\begin{aligned}
T_{\rm state}(t)=\widetilde O\!\Bigg[
&\frac{\mathcal B(t)}{\|P_{\rm tar}\widetilde{\mathbf{Z}}(t)\|}
\Big\{C_{\rm prep}+C_{\rm coeff}+C_{\rm sel}\\
&\quad+C_{\rm be}\big[(\alpha_H+K_{\max}\alpha_L)t+
\log(1/\varepsilon_{\rm sim})\big]\Big\}
\Bigg].
\end{aligned}}
\tag{69}
\]

For $\mathbf b_D=0$, $\delta_{\rm sh}=0$, $P_{\rm tar}=I$, fixed $c_{\rm ker}$ and efficient coefficient access, Eq.~(68a) gives $K_1=O(\log(1/\varepsilon_{\rm alg}))$. Substitution in Eq.~(69) yields the homogeneous query scaling of Ref.~\cite{ref69}, still multiplied by $\lambda_1\|\mathbf Z_0\|/\|\widetilde{\mathbf Z}(t)\|$; it does not prove optimality for an affine lifted PDE. General coefficient loading can cost $O(M_1+M_2R_2)$. For a constant source, the $R_2=2^q$ exponential-bin label can instead be prepared with $q$ rotations, without eliminating source access or controlled simulation. Physical-field preparation additionally requires Eq.~(67a) or \eqref{eq:direct-profile-lcu}.

For a direct coherent implementation of the \(N_t\)-interval recurrence in Algorithm 1, define \(a_{h,j}=e^{\delta_{\rm sh}\Delta t_j}\sum_\ell|\omega^{(1)}_{\ell,j}|\) and let \(a_{s,j}\) be the one-norm of its finite source-time coefficients. One valid LCU normalization obeys

\[
\mathcal B_0=\|\mathbf{Z}_0\|,\qquad
\mathcal B_{j+1}=a_{h,j}\mathcal B_j+a_{s,j}\|\mathbf{b}_D\|,
\qquad
p_{\rm sig,N_t}=
\frac{\|P_{\rm tar}\widetilde{\mathbf{Z}}_{N_t}\|^2}{\mathcal B_{N_t}^2}.
\tag{69a}
\]

Thus its cost contains \(\mathcal B_{N_t}/\|P_{\rm tar}\widetilde{\mathbf{Z}}_{N_t}\|\), not merely \(N_t\max_j a_{h,j}\). If \(a_{h,j}=a>1\), this normalization can grow as \(a^{N_t}\). For the tested kernel weights before shift compensation at \(c_{\rm ker}=1\), \(\varepsilon_{\rm ker}=10^{-8}\), \(K=32\), and \(M=385\), the coefficient one-norm is approximately \(2.373\), so \(a^{16}\) is about \(10^6\). This diagnoses a naive multistep circuit, not a lower bound for all LCHS algorithms. If the target norm tends to zero, no uniform finite bound on this implementation's preparation cost is possible without a signal assumption. One-shot Eq. (69), a normalized positive rule, or a proved stable multistep construction avoids silently assuming linear-in-\(N_t\) success overhead.

Using the saved finite-rule correction states as signal proxies, the normalization-to-signal ratios from Eq. (68) are \(1.15\times10^{18}\) for the one-dimensional constant-\(\mu\) Navier--Stokes order-two run and \(7.81\times10^5\) for the refined two-dimensional compression run. Direct composition of the twelve recorded intervals in Eq. (69a) gives \(1.55\times10^{22}\) and \(1.04\times10^{10}\), respectively. These are normalization-to-archived-signal diagnostics, not measured circuit success probabilities or lower bounds for other implementations. Their size shows that the non-small-signal premise of Theorem~4 remains unverified for these examples despite small classical field errors.

For reconstructed stored NS matrices, outward-rounded weighted bounds certify sufficient shifts \(20,35,11\) for the one-dimensional constant-\(\mu\) orders one/two and refined two-dimensional compression, respectively. These exceed the production shifts and increase their homogeneous normalization factors by about \(3.63\times10^3,1.04\times10^7,9.98\). They are not reruns or certificates of the old shifts; failure of the sufficient test does not establish indefiniteness. A base-operator Lyapunov/block-path alternative remains an unimplemented backend.

If the simulated Hermitian matrices admit efficient row \emph{and} column sparse oracles with at most \(s_{\rm bi}\) entries per row/column and maximum entry magnitude \(a_{\max}\), a standard sparse block encoding has representative normalization and per-call cost

\[
\alpha_H,\alpha_L=O(s_{\rm bi}a_{\max}),
\qquad
C_{\rm be}=\widetilde O(C_{\rm row}+C_{\rm col}+C_{\rm val}+\log D_m).
\tag{70}
\]

The oracle costs on the right must themselves be specified. Equation (51) gives only lifted \emph{row} sparsity and cannot be substituted for \(s_{\rm bi}\). Structured block encodings may avoid column enumeration, but must report their own normalizations and preparation costs.

For the one-shot finite sum, the state and label registers require

\[
\boxed{
S_Q=
O\!\left(
\log D_m+\log(M_1+M_2R_2)
+\log\frac1{\varepsilon_{\rm alg}}+S_{\rm oracle}
\right).
}
\tag{71}
\]

Here \(S_{\rm oracle}\) is the ancilla/workspace cost of matrix, coefficient, source, and observable oracles; a multistep circuit may also need \(O(\log N_t)\) time-label qubits. This qubit count excludes classical memory for oracle data or coefficient tables. Equation (69) prepares an amplitude-encoded state and excludes readout. Let \(T_{\rm state}^{\max}=\max_{t\in\mathcal T_{\rm out}}T_{\rm state}(t)\), where \(\mathcal T_{\rm out}\) contains the \(N_{\rm out}\) requested times, and let \(N_{\rm obs}\) denote the number of bounded observables evaluated at each time. Coherent amplitude estimation repeats controlled state preparation and its inverse, so, suppressing constant and logarithmic overheads,

\[
T_Q^{\rm AE}
=\widetilde O\!\left(
\frac{N_{\rm out}N_{\rm obs}}{\varepsilon_{\rm meas}}
T_{\rm state}^{\max}
\right),
\qquad
T_Q^{\rm samp}
=O\!\left(
\frac{N_{\rm out}N_{\rm obs}}{\varepsilon_{\rm meas}^{2}}
T_{\rm state}^{\max}
\right).
\tag{71a}
\]

Thus measurement cost multiplies, rather than adds to, state-preparation cost. Recovering the physical field at all grid points requires at least \(\Omega(n)\) classical outputs; tomography of the entire lifted state has dimension-dependent sample complexity \cite{ref19}, and dense-output PDE algorithms require separate structure \cite{ref20}. Amplitude estimation requires the preparation circuit and its inverse; physical-field observables additionally require the reconstruction or direct-profile costs specified above.

Signed linear modes require interference, not only an expectation value of the normalized correction. Let $W$ have $N_{\rm obs}$ orthonormal real rows $\mathbf w_i^T$. Prepare the known modal reference with all success ancillas set to zero, orthogonal to the LCU failure branch. For target-preserving scaling, a controlled Hadamard test gives
\[
\begin{aligned}
 p_i&=\tfrac12\left(1+\operatorname{Re}(\mathbf w_i^TP_{\rm tar}\widetilde{\mathbf Z})/\mathcal B\right),\\
 \widehat y_i&=\mathbf w_i^T\bar{\mathbf u}+\mathcal B(2\widehat p_i-1).
\end{aligned}
\tag{71c}\label{eq:signed-modal-readout}
\]
Probability error $\zeta$ adds at most $2\mathcal B\sqrt{N_{\rm obs}}\zeta$; normalized reference-state error $\varepsilon_v$ adds $\mathcal B\sqrt{N_{\rm obs}}\varepsilon_v$. Complex modes need a second phase-shifted test. With a certified external normalization $\overline{\mathcal B}\ge\mathcal B$ and computable $0<s_0\le\|W\mathbf u_{\rm ref}\|$, amplitude estimation~\cite{brassard2002amplitude} has sufficient cost $\widetilde O(N_{\rm obs}^{3/2}\overline{\mathcal B}C_U/(\varepsilon_{\rm obs}s_0))$ for a relative additional backend/readout tolerance; sampling has cost $\widetilde O(N_{\rm obs}^{2}\overline{\mathcal B}^{2}C_U/(\varepsilon_{\rm obs}^{2}s_0^{2}))$. Here $C_U$ includes preparation, simulation, inverse circuits and reflections. This contract does not bound the homotopy/model error. For compressible NS, $\rho=e^q$ is reconstructed classically and is not priced by this linear-output bound.

\subsection{Conditional quantum upper bound and scaling comparison}\label{conditional-quantum-upper-bound-and-scaling-comparison}

For grid-dependent statements, a mesh-consistent choice is
\(\|\mathbf{v}\|_{M_n}^2=\mathbf{v}^*M_n\mathbf{v}\), where a positive mass matrix \(M_n\)
approximates the spatial integral; product blocks use tensor powers of
\(M_n\). On a uniform Cartesian grid, \(M_n=h^{d_x}I\) in spatial
dimension \(d_x\), so a sampled constant has \(O(1)\) norm rather than
the artificial \(O(\sqrt n)\) Euclidean growth. To apply the preceding
Euclidean estimates in this norm, transform each spatial operator by
\(M_n^{1/2}GM_n^{-1/2}\) and consistently transform multilinear maps,
sources, target projectors, and their product-lift blocks. All logarithmic
norms, operator bounds, and block-encoding costs must then be computed
in these coordinates. Mass weighting alone does not uniformly bound
discrete derivatives or homotopy tails; Theorem~4 assumes those estimates,
which the finite-grid experiments do not establish.

For a one-shot call at \(t>0\) with \(\mathbf{b}_D\ne0\), a transparent absolute lifted-state error allocation is
\[
\eta_1=\frac{\varepsilon_{\rm L}}{3\|\mathbf{Z}_0\|},
\qquad
\eta_2=\frac{\varepsilon_{\rm L}}{3t\|\mathbf{b}_D\|},
\qquad
\varepsilon_t=\frac{\varepsilon_{\rm L}}{3\|\mathbf{b}_D\|},
\]
with absent channels omitted and the denominator adjusted accordingly. Choose Eq. (68a) to make the homogeneous and source-integrand defects at most \(\eta_1,\eta_2\), and choose the source-time rule to make its operator defect at most \(\varepsilon_t\). Then the triangle inequality bounds the total finite-rule lifted-state error by \(\varepsilon_{\rm L}\). For a normalized target state with a known lower bound \(\sigma\le\|P_{\rm tar}\mathbf{Z}(t)\|\), the normalized-state error is at most \(2\varepsilon_{\rm L}/\sigma\); bounded-observable error and other coherent errors require a corresponding share of the total budget. The ratios of the input/source norms to \(\sigma\) govern the necessary \(k\)- and time-node precision. Norm and phase estimation for the physical field consume their own budgets. In the multistep emulator, these one-shot allocations are replaced by the amplification-weighted sum.

\textbf{Theorem 4 (conditional sparse-observable scaling).}\\
Fix the polynomial degree \(d\). Suppose, uniformly over the spatial grids under comparison:

\begin{enumerate}
\def\labelenumi{\arabic{enumi}.}
\item
  With a specified mesh-consistent norm on \(X_n\), the homotopy tail has a mesh-independent index \(j_0\), bound \(C_W=\sup_n\|\mathbf{W}_{j_0}^{(n)}\|_{C_T}<\infty\), and ratio \(\|\mathbf{W}_{j+1}^{(n)}\|_{C_T}\le\bar\alpha_{\rm HAM}\|\mathbf{W}_j^{(n)}\|_{C_T}\) for all \(j\ge j_0\), where \(0<\bar\alpha_{\rm HAM}<1\). Its sum solves the same semidiscrete mild equation and lies in the uniqueness ball of Theorem 2. Consequently, \(m=O(\log(1/\varepsilon_{\rm HAM}))\) if \(j_0,C_W,\bar\alpha_{\rm HAM}\) are uniform as specified after Eq. (33).
\item
  The block encodings and their adjoints, normalized \(\mathbf{Z}_0\) and \(\mathbf{b}_D\), Low--Somma coefficients, source-time labels, and target/observable projectors admit coherent access without enumerating \(D_m\) entries. Assume a known lower bound \(0<\sigma\le\|P_{\rm tar}\mathbf{Z}(t)\|\). The costs \(C_{\rm be},C_{\rm prep},C_{\rm coeff},C_{\rm sel},C_{\rm recon},S_{\rm oracle}\) are bounded by \(\operatorname{poly}(m,\log(n+1),1/\varepsilon_{\rm alg})\). The same bound holds for \((\alpha_H+K_{\max}\alpha_L)T\), \(T\|L_{\delta_{\rm sh}}\|\), \(T\Omega_{K_2}\), \(e^{\delta_{\rm sh}T}\), \(\mathcal B(t)/\sigma\), and \(\Gamma_{\rm recon}\).
\item
  At each of \(N_{\rm out}\) requested times, the finite rules and source-time sum meet the allocation displayed above with \(2\varepsilon_{\rm L}/\sigma\) charged to the normalized-state budget, the coherent simulation and reconstruction errors keep the total bounded-observable error below \(\varepsilon_{\rm alg}\), and \(N_{\rm obs}=O(1)\) bounded observables are requested. The preparation circuit and inverse are available for amplitude estimation.
\end{enumerate}

Then the one-shot finite-unitary-sum construction has the conditional bounds

\[
\begin{aligned}
T_Q^{\rm AE}
&=\widetilde O\!\left[
\frac{N_{\rm out}N_{\rm obs}}{\varepsilon_{\rm meas}}
\operatorname{poly}\!\left(
m,d,\log(n+1),\frac1{\varepsilon_{\rm alg}}
\right)\right],\\
S_Q
&=\operatorname{poly}\!\left(
m,d,\log(n+1),\log\frac1{\varepsilon_{\rm alg}}
\right)+S_{\rm oracle}.
\end{aligned}
\tag{71b}
\]

By contrast, the two separate tasks of writing and materializing the \emph{full lifted state} satisfy
\[
T_C^{\rm full\ output}=\Omega(D_m),\qquad
S_C^{\rm materialized\ state}=\Omega(D_m).
\tag{72}
\]

\textbf{Proof.} The geometric tail bound following Eq. (33) gives the stated \(m\). Equation (50c) implies \(\log D_m=O_d(m\log(n+1))\). Under item 2, Eq. (68a) gives logarithmic cutoff and at most polynomially many \(k\)- and source-time labels in \(1/\varepsilon_{\rm alg}\); their \emph{coherent} preparation cost is separately bounded in that item. Equations (69) and (71) then give the claimed dimension and workspace dependence, and Eq. (71a) supplies the multiplicative measurement factor. The two bounds in Eq. (72) follow solely from the number of written or stored scalars. \(\square\)

Theorem 4 gives a conditional sparse-observable upper bound. Equation (72) instead concerns full lifted-state output, using \(D_{\rm sym}\) when both sides apply exact compression. Assessing a same-observable quantum advantage requires comparison with direct classical PDE solvers and accounting for profile generation, data loading, source-time quadrature and output preparation. The current theorem supplies the quantum-side conditions, not a matched classical query lower bound.

\subsection{Dependence on spatial resolution}\label{dependence-on-spatial-resolution}

For a one-dimensional grid with \(n\) points and \(h=O(n^{-1})\), the derivative, frozen-operator and unscaled lift bounds are

\begin{align}
\|D_1\|&=O(n),\quad \|D_2\|=O(n^2),\quad \|D_3\|=O(n^3), \tag{73}\\
\|G_B\|&=O\!\left(\nu n^2+\|\bar u\|_\infty n+\|\bar u_x\|_\infty\right), \tag{74}\\
\|G_K\|&=O\!\left(n^3+(|c|+6\|\bar u\|_\infty)n+6\|\bar u_x\|_\infty\right), \tag{75}\\
\|A_m\|&\le p_{\max}\|G\|+\|A_{\rm coupling}\|. \tag{76}
\end{align}

These bounds are not block-encoding constructions. For standard finite differences, \(\|D_2\|=\Theta(n^2)\) and a nondegenerate third-derivative stencil has \(\|D_3\|=\Theta(n^3)\); therefore the sufficient Hamiltonian-simulation cost in Eq. (69) is not polylogarithmic in grid resolution merely because its state register is. Similarity scaling, spectral shifting, and the chosen encoding can further increase normalization or reduce the success amplitude. For example, the second-order constant-dynamic-viscosity one-dimensional Navier--Stokes run uses \(\delta_{\rm sh}\simeq24.23\) over \(T=1.5\), giving \(e^{\delta_{\rm sh}T}\simeq6.1\times10^{15}\) in Eq. (68); this is a normalization factor, not a measured failure probability. Theorem 4's access and signal assumptions are thus not verified by the reported finite-rule field errors. Preconditioning, spectral bases, an interaction picture, or structured block encodings might improve the dependence, but each requires a separate construction and analysis.

% Local numerical-section layout; manuscript wording and figure data are unchanged.
\begingroup
\setlength{\tabcolsep}{4pt}
\renewcommand{\arraystretch}{1.12}
\setlength{\abovecaptionskip}{7pt}
\setlength{\belowcaptionskip}{5pt}
\renewcommand{\topfraction}{0.95}
\renewcommand{\textfraction}{0.06}
\renewcommand{\floatpagefraction}{0.75}
\makeatletter
\setlength{\@fptop}{0pt}
\setlength{\@fpsep}{12pt}
\setlength{\@fpbot}{0pt plus 1fil}
\makeatother
\section{Numerical experiments}\label{numerical-experiments}

All plotted Frozen--Oseen fields are classical emulations of the equidistant Low--Somma rule with $K_1=K_2=32$, $M_1=M_2=385$, $c=1$, and $\varepsilon_{\rm ker}=10^{-8}$; the profile remains fixed. The active symmetric-lift dimension $D_{\rm sym}$ is padded to $N_Q=2^{\lceil\log_2D_{\rm sym}\rceil}$, including one affine-source coordinate. This counts state-register qubits only, not ancillary qubits or a quantum execution.

The complete-lift diagnostic is
\[
E_{\rm lift}=\frac{\|\mathbf Y_K-\mathbf Y_{\rm dir}\|_2}{\|\mathbf Y_{\rm dir}\|_2},
\qquad \mathbf Y_K=D\mathbf Z_K,\quad \mathbf Y_{\rm dir}=D\mathbf Z_{\rm dir}.
\tag{76a}\label{eq:complete-lift-error}
\]
Both states use the same compressed coordinates after undoing tensor scaling; physical-field errors have separately stated denominators.

\paragraph{Comparison protocol}\label{fairprotocol}
Each comparison fixes the physical target, grid, initial data, and output norm. Conventional HAM and IQHAM are classically integrated; matched auxiliary controls additionally fix the zeroth coefficient and $\hbar$.

\subsection{Forced Burgers equation}\label{forced-burgers-equation-analysis-and-numerical-diagnostics}

\subsubsection{Governing equation and exact Oseen split}\label{governing-equation-and-exact-oseen-split}

We consider the classical viscous Burgers equation \cite{ref21,ref22}, here augmented by a stationary body force,

\[
u_t+uu_x=\nu u_{xx}+\sin(\pi x),
\qquad 0<x<1,
\tag{77}
\]

\[
u(0,t)=u(1,t)=0,\qquad
u(x,0)=0.3\sin(\pi x),
\tag{78}
\]

with

\[
\nu=0.01,\qquad T=1.
\tag{79}
\]

For

\[
\mathcal F_B(u)=\nu u_{xx}-uu_x,
\qquad f_B(x)=\sin(\pi x),
\tag{80}
\]

the full Fréchet--Oseen operator at \(\bar u\) is

\[
\boxed{
\mathcal G_Bv
=\nu v_{xx}-\bar u\,v_x-\bar u_xv.
}
\tag{81}
\]

The residual and exact nonlinear remainder are

\[
r_B
=\nu\bar u_{xx}-\bar u\bar u_x+\sin(\pi x),
\qquad
\mathcal R_B(w)=-ww_x.
\tag{82}
\]

Thus

\[
w_t=\mathcal G_Bw+r_B-ww_x.
\tag{83}
\]

With \(\mathscr L_B=\partial_t-\mathcal G_B\), the homotopy equations are

\[
\mathscr L_BW_0=r_B,\qquad
W_0(x,0)=0.3\sin(\pi x)-\bar u(x),
\tag{84}
\]

\[
\boxed{
\mathscr L_BW_j
=-\sum_{a+b=j-1}W_a\,\partial_xW_b,
\qquad
W_j(x,0)=0,\quad j\ge1.
}
\tag{85}
\]

For homogeneous Dirichlet data, integration by parts and the Poincaré inequality give

\[
\mu_2(\mathcal G_B)
\le-\nu\pi^2
+\frac12\|(\bar u_x)_-\|_\infty,
\tag{87}
\]

where \(a_-=\max(-a,0)\). A sharp negative profile slope can reduce the residual while increasing linear amplification, so the fitted layer width and propagator must be checked together \cite{ref23}.

\subsubsection{Frozen profile}\label{formal-asymptotic-derivation-of-the-right-boundary-layer-profile}
\label{frozen-profile-used-in-the-proof-of-concept-code}

For the steady balance \(\nu\bar u''-\bar u\bar u'+\sin(\pi x)=0\), the outer approximation is
\[
\bar u_{\rm out}(x)=\frac{2}{\sqrt\pi}\sin\frac{\pi x}{2},
\qquad U_R=\bar u_{\rm out}(1)=\frac{2}{\sqrt\pi}.
\tag{88}
\]
The right-layer coordinate \(y=(1-x)/\nu\) gives \(U_{yy}+UU_y=0\), $U(0)=0$, $U(\infty)=U_R$, and
\[
U(y)=U_R\tanh\frac{U_Ry}{2}.
\tag{89}
\]
This balance motivates the $O(\nu)$ width and tanh factor, but not the fitted constants \cite{ref41,ref44}.

The implemented profile uses a multiplicative boundary-layer enrichment:

\[
\boxed{
\bar u_\theta(x)
=s\frac{2}{\sqrt\pi}
\sin\frac{\pi x}{2}
\tanh\!\left(\frac{1-x}{\delta_{\rm BL}}\right)
+a\sin(\pi x),
}
\tag{90}
\]

\[
\delta_{\rm BL}=b\frac{2\nu}{s(2/\sqrt\pi)}.
\tag{91}
\]

The stored run used

\[
s=0.905,\qquad b=1.064,\qquad a=0.090,
\qquad \delta_{\rm BL}=2.083857345\times10^{-2}.
\tag{92}
\]

The profile satisfies both boundaries. Its parameters were fitted once on this benchmark and held fixed across orders; the main-case errors are therefore in-sample.

\subsubsection{Discretization and results}\label{discretization-and-numerical-baseline}

The lifted test has $n=5$ interior correction points and $\hbar=-1$. Diffusion uses centered differences; transport uses backward differences. Profile derivatives are evaluated analytically from Eq.~(90), giving

\[
\dot{\mathbf w}=
\nu D_2\mathbf w-\mathbf p\odot D_1\mathbf w-\mathbf p_x\odot\mathbf w
+\nu\mathbf p_{xx}-\mathbf p\odot\mathbf p_x+\mathbf f-\mathbf w\odot D_1\mathbf w,
\tag{92a}
\]

Here $D_1,D_2$ are difference matrices and $\odot$ is componentwise multiplication. The analytic profile derivatives make this an enriched semidiscrete model; $n=5$ does not resolve the $O(\nu)$ layer.

Equivalently, the reconstructed nodal variable $\mathbf u=\mathbf p+\mathbf w$ solves
\[
\begin{aligned}
\dot{\mathbf u}=\mathcal F_{h,\mathbf p}(\mathbf u)
={}&\nu D_2\mathbf u-\mathbf u\odot D_1\mathbf u+\mathbf f\\
&+\nu(\mathbf p_{xx}-D_2\mathbf p)
+\mathbf u\odot(D_1\mathbf p-\mathbf p_x).
\end{aligned}
\tag{92b}\label{eq:enriched-burgers}
\]
The frozen operator is the Jacobian of this fixed enriched model at $\mathbf p$; varying the enrichment would change the spatial model. All curves reconstruct the analytic profile plus a piecewise-linear correction. The same-grid nonlinear reference differs from a 400-node PDE proxy by $8.51\times10^{-3}$.

All three methods use the same five-point semidiscrete target and nonlinear nodal reference, with error
\[
e_{\rm node}=\frac{\|\mathbf u_h-\mathbf u_{h,\rm ref}\|_2}{\|\mathbf u_{h,\rm ref}\|_2}.
\tag{93}
\]
The Frozen--Oseen hierarchy uses 16 fixed linear substeps and exact symmetric-product compression. The other two hierarchies are integrated classically; their curves are not LCHS results.

\begin{table}[!htbp]
\centering
\caption{Common-grid Burgers comparison at $T=1$. $D_{\rm sym}$ is the active lifted dimension, $N_Q$ its exact register embedding, and both errors are relative to the same $n=5$ nonlinear semidiscrete reference. The conventional HAM uses a diffusion auxiliary and classical hierarchy propagation.}
\label{tab:burgers-common-grid}
\TableFont
\begin{tabular*}{\linewidth}{@{\extracolsep{\fill}}rrrrr@{}}
\toprule
$m$ & $D_{\rm sym}$ & $N_Q$ & Frozen--Oseen LCHS & Conventional HAM\\
\midrule
1 & 30 & 32 & $1.05269\times10^{-2}$ & $3.096\times10^{-1}$\\
2 & 95 & 128 & $4.87970\times10^{-3}$ & $2.033\times10^{-1}$\\
3 & 285 & 512 & $2.37859\times10^{-3}$ & $1.349\times10^{-1}$\\
4 & 791 & 1,024 & $1.20068\times10^{-3}$ & $9.453\times10^{-2}$\\
5 & 2,056 & 4,096 & $6.19361\times10^{-4}$ & $7.388\times10^{-2}$\\
\bottomrule
\end{tabular*}
\end{table}

At order five, Frozen--Oseen HAM-LCHS has nodal error $6.19\times10^{-4}$ versus $7.39\times10^{-2}$ for diffusion-auxiliary HAM (Table~\ref{tab:burgers-common-grid}). IQHAM~\cite{ref12} reaches $5.69\times10^{-4}$ after 26 outer updates: it improves the diffusion auxiliary by iteration, whereas the present method changes the fixed auxiliary. The complete-lift LCHS/direct discrepancy is $1.13\times10^{-6}$.

\begin{figure}[!htbp]
\centering
\includegraphics[width=\linewidth]{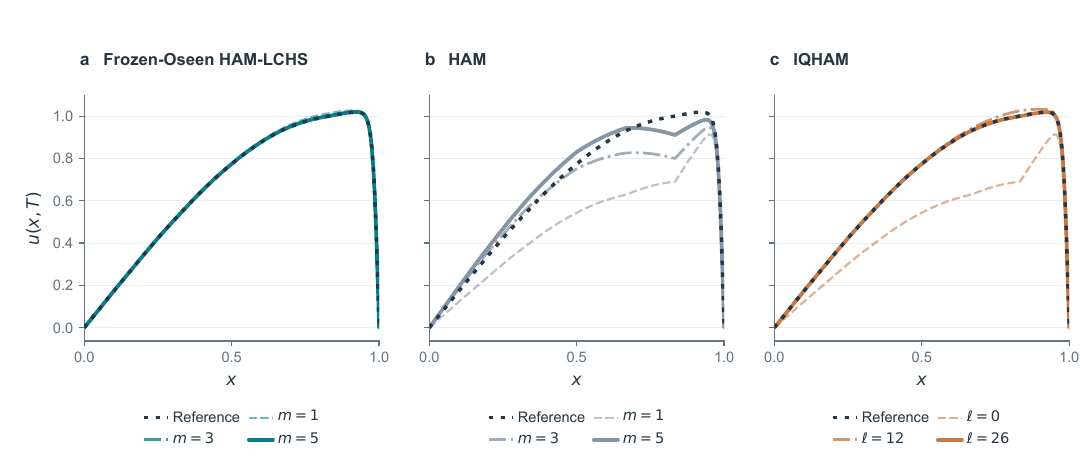}
\caption{Burgers profiles at $T=1$ on the common $n=5$ semidiscrete grid. (a) Finite-$K_1/K_2$ Frozen--Oseen HAM-LCHS at orders $m=1,3,5$. (b) Conventional diffusion-auxiliary HAM, integrated classically at the same orders. (c) IQHAM denotes the iterative QHAM introduced by Xue et al.~\cite{ref12}; its outer indices $\ell=0,12,26$ are integrated classically here. The dotted line is the same-grid nonlinear reference. The $n=5$ grid does not resolve the boundary layer.}
\label{fig:burgers-three-methods}
\end{figure}

\paragraph{Validation}\label{burgers-error-diagnostics}

Four unfitted $n=17$ holdouts with $\bar u=u_{\rm in}$ improve at order two but remain grid-underresolved: their $n=17$ to $129$ discrepancies are $5.23\times10^{-2}$--$8.53\times10^{-2}$.

\Needspace{5\baselineskip}
\paragraph{Independent auxiliary-operator ablation}
For the four unfitted parameter pairs, we compare diffusion, mean transport $\nu D_2-\langle\mathbf p\rangle D_1$, frozen advection $\nu D_2-\operatorname{diag}(\mathbf p)D_1$, and the full Jacobian $J=DF_n(\mathbf p)$ with $\mathbf p=\mathbf u_{\rm in}$. The same centered grids, $\hbar=-1$, and initial coefficient are used. Table~\ref{tab:auxiliary-ablation} reports classical hierarchy diagnostics, not additional LCHS curves.

\begin{table}[!htbp]
\centering\TableFont
\caption{Unfitted Burgers auxiliary-operator ablation at $N_x=127$, $m=5$, $T=0.75$. Entries are relative nodal errors against the same-grid nonlinear reference. The last column compares that reference with a 511-node nonlinear proxy. All rows are retained.}
\label{tab:auxiliary-ablation}
\begin{tabular*}{\linewidth}{@{\extracolsep{\fill}}lrrrrr@{}}
\toprule
$(a,f)$ & Diffusion & Mean transport & Frozen advection & Full $J$ & Spatial proxy\\
\midrule
$(0.2,0.8)$ & $4.107\,10^{-2}$ & $1.124\,10^{-2}$ & $7.643\,10^{-3}$ & $4.550\,10^{-3}$ & $1.141\,10^{-4}$\\
$(0.2,1.2)$ & $2.137\,10^{-1}$ & $8.366\,10^{-2}$ & $6.397\,10^{-2}$ & $5.055\,10^{-2}$ & $3.587\,10^{-4}$\\
$(0.4,0.8)$ & $2.830\,10^{-1}$ & $2.775\,10^{-2}$ & $2.015\,10^{-2}$ & $2.227\,10^{-2}$ & $5.187\,10^{-4}$\\
$(0.4,1.2)$ & $9.157\,10^{-1}$ & $1.600\,10^{-1}$ & $1.172\,10^{-1}$ & $1.770\,10^{-1}$ & $1.090\,10^{-3}$\\
\bottomrule
\end{tabular*}
\end{table}

Transport improves all four fifth-order controls over diffusion, but full $J$ is not uniformly best (Table~\ref{tab:auxiliary-ablation}). In the first case it reaches $1\%$ error one order earlier than frozen advection, reducing the symmetric-lift dimension by factors $7.46$, $12.94$, and $23.71$ on the three tested grids; the other cases do not reach that target through order five. These are representation counts from classical hierarchy diagnostics, not measured quantum speedups. The short-time sufficient test does not certify the full $T=0.75$ interval.
\FloatBarrier

\subsection{Korteweg--de Vries equation}\label{kdv-profile-selection-and-numerical-diagnostics}

\subsubsection{Periodic traveling-wave benchmark}
We consider
\[
u_t+6uu_x+u_{xxx}=0,\qquad x\in[-\pi,\pi],\qquad 0\le t\le1,
\tag{94}
\]
with periodic boundary conditions. An exact cnoidal wave is
\[
u_*(x,t)=a+2\eta\beta^2\operatorname{cn}^2\bigl(\beta(x-ct);\eta\bigr),
\quad \beta=\frac{\mathrm K(\eta)}{\pi},
\quad c=6a+4\beta^2(2\eta-1),
\tag{95}
\]
Here $\eta=0.1$ is the elliptic parameter, $\mathrm K$ the complete elliptic integral, and $\operatorname{cn}$ the Jacobi elliptic cosine. We set $a=0.5$ and $u_{\rm in}=u_*(\cdot,0)$, giving $c=2.1570184182$ \cite{bottman2009cnoidal}; the computation stays in laboratory coordinates.

Set \(\bar u=u_{\rm in}\) once and write \(u=\bar u+w\). The exact correction equation is
\[
w_t=\mathcal G_Kw+r_K-6ww_x,\qquad w(x,0)=0,
\tag{96}
\]
where
\[
\mathcal G_Kw=-w_{xxx}-6\bar u w_x-6\bar u_xw,
\qquad r_K=-6\bar u\bar u_x-\bar u_{xxx}=-c\bar u_x.
\tag{97}
\]
The initial-state profile is fixed and has a nonzero residual, so the $K_2$ source channel is used. With $\mathscr L_K=\partial_t-\mathcal G_K$, the hierarchy is
\[
\mathscr L_KW_0=r_K,\qquad
\mathscr L_KW_j=-6\sum_{p+q=j-1}W_p\partial_xW_q,
\quad W_j(x,0)=0\quad(j\ge0),
\tag{98}
\]
and \(u^{[m]}=\bar u+\sum_{j=0}^mW_j\).

Although the wave translates, its small modulation makes the initial profile useful: if $b$ is its mean, translation invariance gives
\[
\|u_*(\cdot,t)-\bar u\|_{H^s}
\le2\|\bar u-b\|_{H^s}=O(\eta),\qquad \eta\to0.
\tag{99}
\]
Thus the frozen operator retains mean transport and variable-profile coupling while the remainder is $O(\eta^2)$ for this family; a dispersion-only auxiliary leaves first-order transport in the forcing. This local comparison does not establish uniform convergence. On a fixed-grid ball of radius $\rho$,
\[
\|{-6\mathbf w\odot D_1\mathbf w}+6\mathbf z\odot D_1\mathbf z\|_2
\le12\rho\|D_1\|_2\|\mathbf w-\mathbf z\|_2.
\tag{100}
\]
The convergence condition also depends on the frozen propagator (Section~2). For smooth periodic functions,
\[
\operatorname{Re}\langle w,\mathcal G_Kw\rangle
=-3\int_{-\pi}^{\pi}\bar u_x|w|^2\,dx
\le3\|\bar u_x\|_\infty\|w\|_2^2.
\tag{101}
\]
Because $\bar u_x=O(\eta)$, this bounds finite-time growth but does not imply contractivity of the laboratory-frame operator \cite{bottman2009cnoidal}.

\subsubsection{Discretization and error control}
We use seven-point periodic Fourier collocation with pointwise products. The residual and frozen operator use the same discrete map,
\[
F_n(\mathbf v)=-D_3\mathbf v-6\mathbf v\odot D_1\mathbf v,\qquad
\mathbf r_n=F_n(\bar{\mathbf u}_n),\qquad G_n=DF_n(\bar{\mathbf u}_n).
\tag{102}
\]
The correction is Fourier-interpolated onto 512 comparison points and added to the analytic initial profile, without phase alignment. The wave-normalized metric is
\[
e_{\rm wave}=\frac{\|\mathbf u_{\rm num}-\mathbf u_*\|_2}{\|\mathbf u_*-b_{512}\mathbf 1\|_2}.
\tag{103}
\]
Here $b_{512}$ is the sampled exact-wave mean. The seven-point nonlinear spatial error is $1.5663\times10^{-5}$; nine and 13 points reduce it to $6.3359\times10^{-8}$ and $3.9201\times10^{-12}$ for this smooth wave.

The separate propagation metric is
\[
E_{\rm prop,wave}=\frac{\|\mathbf u_{\rm LCHS}-\mathbf u_{\rm dir}\|_2}
{\|\mathbf u_*-b_{512}\mathbf1\|_2}.
\]
This real-field discrepancy is distinct from the complete complex-lift error $E_{\rm lift}$.

The exact symmetric-product compression retains the complete hierarchy. Finite LCHS uses four fixed intervals and tensor scale $0.06$.

\begin{table}[!htbp]
\centering
\caption{Cnoidal-wave lift sizes and errors. $N_Q$ is the exact zero-padded register dimension. Every order uses \(K_1=K_2=32\), 385 equidistant rectangle nodes, and four intervals. Times are single-run propagation costs of a three-worker classical emulator with one BLAS thread per worker, excluding construction and direct diagnostics. The last column is $E_{\rm prop,wave}$ for real-projected physical fields reconstructed on 512 points, not $E_{\rm lift}$.}
\label{tab:kdv-finite-k}
\TableFont
\textbf{(a) Lift dimensions and propagation cost}\par\smallskip
\begin{tabular*}{\linewidth}{@{\extracolsep{\fill}}rrrrrr@{}}
\toprule
\(m\) & Full dimension & Symmetric dimension & $N_Q$ & Symmetric nonzeros & Time (s) \\
\midrule
1 & 70 & 49 & 64 & 602 & 2.243 \\
2 & 518 & 189 & 256 & 3,255 & 4.953 \\
3 & 4,102 & 679 & 1,024 & 14,721 & 14.353 \\
4 & 32,774 & 2,226 & 4,096 & 58,002 & 27.409 \\
5 & 262,150 & 6,748 & 8,192 & 204,936 & 87.094 \\
\bottomrule
\end{tabular*}
\par\smallskip
\textbf{(b) Wave-normalized errors}\par\smallskip
\begin{tabular*}{\linewidth}{@{\extracolsep{\fill}}rccc@{}}
\toprule
\(m\) & Direct/exact & LCHS/exact & $E_{\rm prop,wave}$ \\
\midrule
1 & \(2.843\times10^{-3}\) & \(2.843\times10^{-3}\) & \(6.164\times10^{-9}\) \\
2 & \(1.587\times10^{-4}\) & \(1.587\times10^{-4}\) & \(4.424\times10^{-9}\) \\
3 & \(2.287\times10^{-5}\) & \(2.287\times10^{-5}\) & \(3.469\times10^{-9}\) \\
4 & \(1.551\times10^{-5}\) & \(1.551\times10^{-5}\) & \(5.195\times10^{-9}\) \\
5 & \(1.565\times10^{-5}\) & \(1.565\times10^{-5}\) & \(1.744\times10^{-9}\) \\
\bottomrule
\end{tabular*}
\end{table}

Wave error reaches the $1.57\times10^{-5}$ spatial-error level by order four, and all LCHS/direct wave discrepancies are below $6.2\times10^{-9}$ (Table~\ref{tab:kdv-finite-k}).

\Needspace{5\baselineskip}
\subsubsection{Matched auxiliary controls and HAM recipes}
At common $\hbar=-1$ and zeroth-coefficient initialization, matched classical hierarchy controls isolate auxiliary choice (Table~\ref{tab:kdv-matched-auxiliary}). Mean transport uses $A=-D_3-6bD_1$, where $b$ is the known initial-wave mean; frozen advection uses $A=-D_3-6\operatorname{diag}(\bar{\mathbf u}_n)D_1$. Mean transport reaches the $10^{-4}$ PDE-wave target at order two, one order before full $G_n$; both approach the spatial-error level by order five. These controls do not compare quantum runtimes.

\begin{table}[!tp]
\centering\TableFont
\caption{Matched KdV auxiliary control: PDE-wave errors on the same seven-point system, forced zeroth coefficient and $\hbar=-1$. $m_*$ is the first tested order reaching $10^{-4}$; a dash means none through order five. All entries are classical coefficient diagnostics.}
\label{tab:kdv-matched-auxiliary}
\begin{tabular*}{\linewidth}{@{\extracolsep{\fill}}lcccc@{}}
\toprule
Auxiliary & $m=2$ & $m=3$ & $m=5$ & $m_*$\\
\midrule
Dispersion & $2.512$ & $1.699$ & $8.492\,10^{-1}$ & --\\
Mean transport & $7.822\,10^{-5}$ & $1.895\,10^{-5}$ & $1.568\,10^{-5}$ & 2\\
Frozen advection & $3.772\,10^{-4}$ & $3.088\,10^{-5}$ & $1.571\,10^{-5}$ & 3\\
Full Jacobian & $1.587\,10^{-4}$ & $2.287\,10^{-5}$ & $1.565\,10^{-5}$ & 3\\
\bottomrule
\end{tabular*}
\end{table}

In the separate recipe comparison (Fig.~\ref{fig:kdv-three-methods}), dispersion-auxiliary HAM has order-five wave error $0.66028$; IQHAM~\cite{ref12} reaches $1.5626\times10^{-5}$ after nine outer updates. These recipes also differ in the zeroth trajectory and $\hbar$; the matched controls above isolate only the auxiliary.

\begin{figure}[!htbp]
\centering
\includegraphics[width=\linewidth]{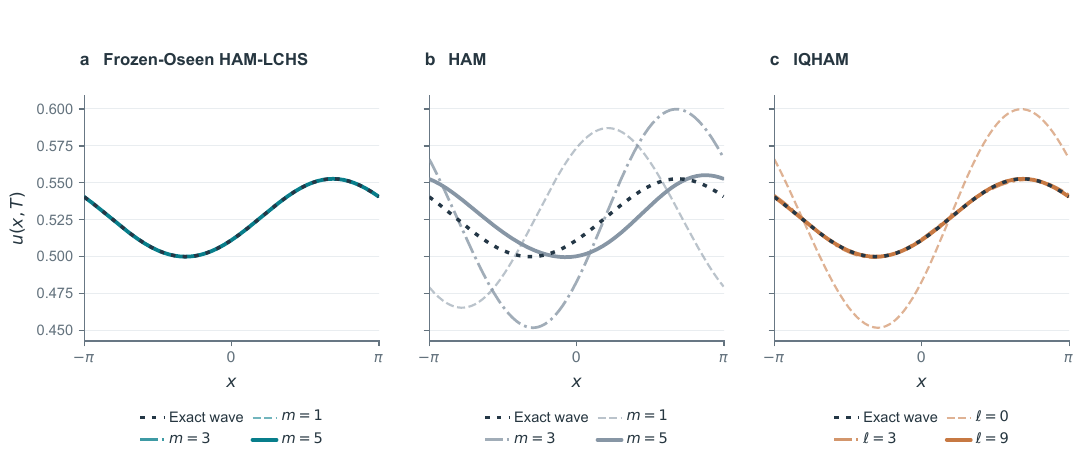}
\caption{KdV recipe comparison at $T=1$ on the common $n=7$ Fourier grid. (a) Finite-$K_1/K_2$ Frozen--Oseen HAM-LCHS at $m=1,3,5$. (b) Conventional dispersion-auxiliary HAM at the same orders. (c) IQHAM, the iterative QHAM of Xue et al.~\cite{ref12}, at $\ell=0,3,9$. Panels (b,c) use classical hierarchy integration; their zeroth trajectory and $\hbar$ differ from (a). The dotted line is the exact cnoidal wave, with common axes and no phase alignment. Table~\ref{tab:kdv-matched-auxiliary} isolates the auxiliary choice.}
\label{fig:kdv-three-methods}
\end{figure}
\FloatBarrier

\paragraph{Representation and resource tradeoff}
At the common $10^{-4}$ target, mean transport uses a smaller complete lift and lower output-aware query proxy than full $G_n$ under the stated access model (Table~\ref{tab:kdv-resource-tradeoff}). The proxy is not a gate count or measured runtime.

\begin{table}[!htbp]
\centering\TableFont
\caption{KdV resource diagnostic at a common $10^{-4}$ PDE-wave threshold, not identical achieved errors. Both systems use scale $0.06$, $K_1=K_2=48.3471$, 1077 raw rectangle nodes and one interval. $\mathcal B$ is the source-branch normalization. $Q_{\rm proxy}$ concerns two signed Fourier modes at the same additional readout tolerance; it is not a gate count or runtime.}
\label{tab:kdv-resource-tradeoff}
\begin{tabular*}{\linewidth}{@{\extracolsep{\fill}}lrrrccc@{}}
\toprule
Auxiliary & $m$ & $D_{\rm sym}$ & $N_Q$ & PDE error & $\mathcal B$ & $Q_{\rm proxy}$\\
\midrule
Mean transport & 2 & 189 & 256 & $7.822\,10^{-5}$ & $40.309$ & $1.048\,10^{13}$\\
Full Jacobian & 3 & 679 & 1024 & $2.287\,10^{-5}$ & $311.384$ & $1.298\,10^{14}$\\
\bottomrule
\end{tabular*}
\end{table}

\FloatBarrier
\subsection{Two-dimensional Burgers and Zakharov--Kuznetsov equations}
\label{sec:two-dimensional-tests}

On the periodic square $[-\pi,\pi]^2$, we consider the vector Burgers equation
\[
 \boldsymbol U_t+(\boldsymbol U\cdot\nabla)\boldsymbol U=0.1\Delta\boldsymbol U,\qquad \boldsymbol U=(u,v)^T,
 \tag{104}\label{eq:burgers2d-test}
\]
with the smooth exact solution $\boldsymbol U=\boldsymbol b-0.2\nabla\log\phi$, where
\[
 \begin{split}
 \boldsymbol b&=(0.8,0.4)^T,\qquad \xi=x-0.8t,\quad \eta=y-0.4t,\\
 \phi&=1+0.15e^{-0.1t}(\cos\xi+\cos\eta)
       +0.05e^{-0.2t}\cos\xi\cos\eta.
 \end{split}
 \tag{105}
\]
This is a smooth Cole--Hopf solution with $\phi\ge0.65$; we use its initial data and $T=1$ \cite{ref67}. The second problem is the Zakharov--Kuznetsov (ZK) equation~\cite{ref68},
\[
 \begin{aligned}
 u_t+6uu_x+\partial_x\Delta u&=0,\\
 u(x,y,0)&=0.5+0.05\cos x+0.04\cos y+0.03\cos(x+y).
 \end{aligned}
 \tag{106}
 \label{eq:zk-test}
\]
at $T=0.5$. Both profiles are frozen at the initial state. The operators are $\mathcal G\boldsymbol Z=0.1\Delta\boldsymbol Z-(\bar{\boldsymbol U}\cdot\nabla)\boldsymbol Z-(\boldsymbol Z\cdot\nabla)\bar{\boldsymbol U}$ for Burgers and $\mathcal Gz=-\partial_x\Delta z-6\bar u z_x-6\bar u_xz$ for ZK.

\Needspace{5\baselineskip}
Both tests use $9\times9$ Fourier grids, first-order lifts, four fixed linear substeps, and tensor scale $0.1$. Burgers uses the exact solution as reference; ZK uses a refined nonlinear reference. For these and the Navier--Stokes (NS) tests, we report
\[
 E_{\rm fluc}=
 \frac{\|\mathbf z_h-\mathbf z_{\rm ref}\|_2}
 {\|\mathbf z_{\rm ref}-\langle\mathbf z_{\rm ref}\rangle\|_2},
 \tag{107}
 \label{eq:extension-fluctuation-error}
\]
where the numerator stacks the sampled physical-field differences and each reference component is centered only in the denominator. First-order $E_{\rm fluc}$ is $3.49\times10^{-5}$ for Burgers and $7.03\times10^{-4}$ for ZK (Table~\ref{tab:extension-errors}); the corresponding complete-lift LCHS/direct discrepancies are $3.06\times10^{-9}$ and $4.62\times10^{-9}$. Figure~\ref{fig:two-dimensional-lchs} shows the full fields and signed errors.

\begin{figure}[p]
\centering
\includegraphics[width=\linewidth]{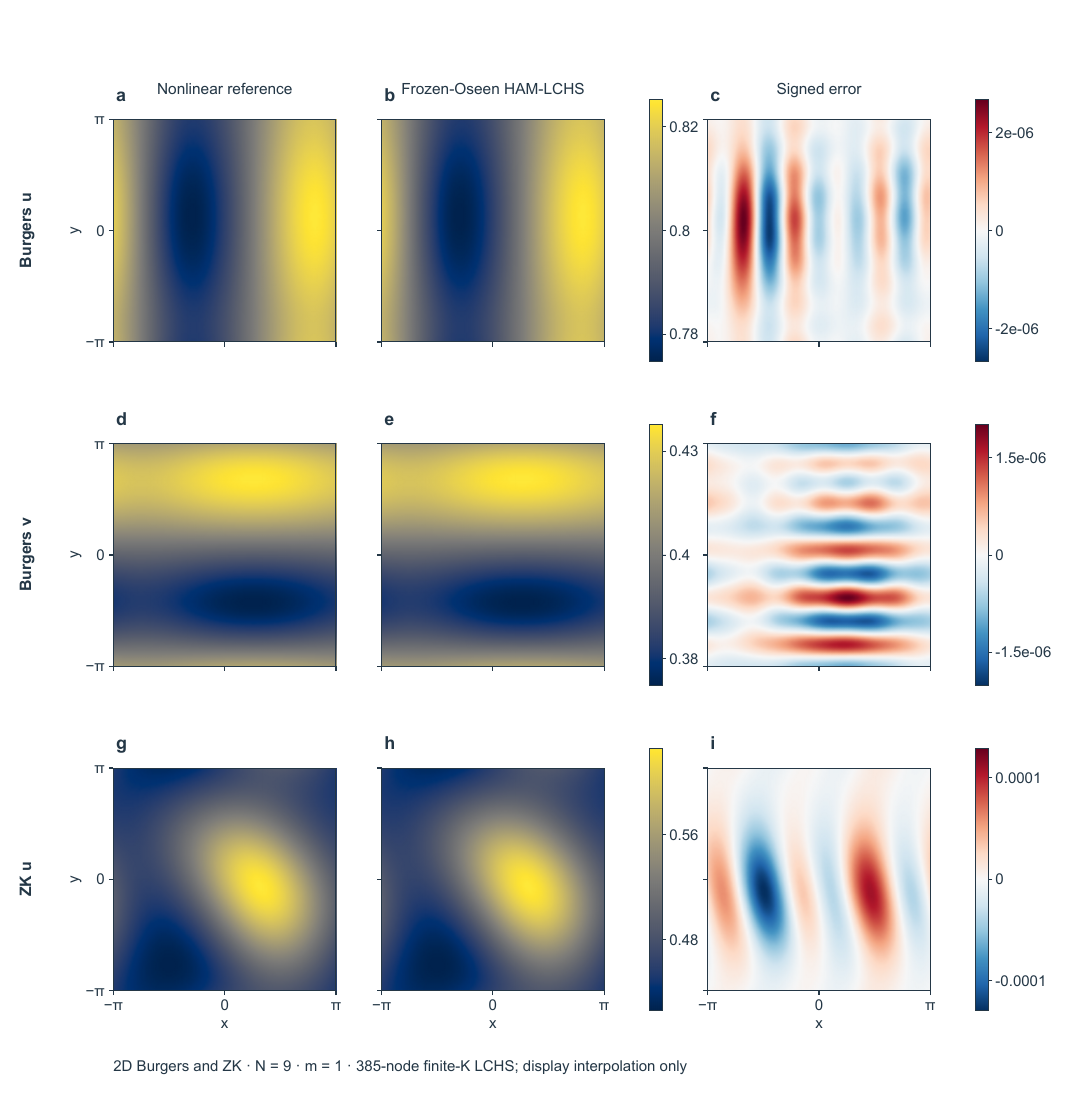}
\caption{Two-dimensional Burgers and ZK fields at $T=1$ and $T=0.5$, respectively, for Eqs.~\eqref{eq:burgers2d-test} and~\eqref{eq:zk-test}. Rows show Burgers $u$, Burgers $v$, and ZK $u$; columns show the reference, first-order finite-$K_1/K_2$ Frozen--Oseen HAM-LCHS result, and pointwise signed error $e_\phi=\phi_{\rm LCHS}-\phi_{\rm ref}$ for the displayed field $\phi$. Each reference/result pair shares a scale; error maps use separate zero-centered scales. The aggregate error $E_{\rm fluc}$ is defined in Eq.~\eqref{eq:extension-fluctuation-error}. Resampling to $128\times128$ is for display only.}
\label{fig:two-dimensional-lchs}
\end{figure}

\subsection{One-dimensional isothermal compressible Navier--Stokes equations}
\label{sec:ns-tests}

We solve
\[
 \rho_t+(\rho u)_x=0,\qquad
 (\rho u)_t+(\rho u^2+\rho)_x=(\eta(\rho)u_x)_x
 \tag{108}
 \label{eq:ns-test}
\]
on $[-\pi,\pi]$ with periodic boundaries, $T=1.5$,
$\rho(x,0)=e^{0.3\cos x}$, and $u(x,0)=-0.6\sin x$.
We distinguish constant kinematic viscosity, $\eta(\rho)=\nu\rho$,
from constant dynamic viscosity, $\eta(\rho)=\mu$, with
$\nu=\mu=0.005$. With $q=\log\rho$, the former is quadratic:
\[
 q_t=-u q_x-u_x,\qquad
 u_t=-u u_x-q_x+\nu(u_{xx}+q_xu_x).
 \tag{109}
\]
For constant $\mu$, introducing $s=e^{-q}$ instead gives the polynomial system
\[
 \begin{aligned}
 q_t&=-u q_x-u_x, & s_t&=suq_x+su_x,\\
 u_t&=-u u_x-q_x+\mu s u_{xx}.
 \end{aligned}
 \tag{110}
\]
Its cubic terms are included in the product lift.

The fixed profile is obtained from the linear acoustic--viscous predictor
\[
 q^a_t=-u^a_x,\qquad u^a_t=-q^a_x+0.005u^a_{xx},
 \qquad (q^a,u^a)(0)=(q,u)(0).
 \tag{111}
\]
Writing $\boldsymbol V^a=(q^a,u^a)^T$, we use the time average
$(\bar q,\bar u)^T=T^{-1}\int_0^T\boldsymbol V^a(t)\,dt$, evaluated by 16-point Gauss quadrature.
For constant $\mu$, we also test $(\bar q,\bar u)^T=[\boldsymbol V^a(0)+\boldsymbol V^a(T)]/2$;
in both cases $\bar s=e^{-\bar q}$ is formed afterwards.
The time average is the constant profile closest in integrated squared distance to the linear predictor; it uses no nonlinear reference. Because it differs from the initial state, the complete lifted initial products include $\mathbf W_0(0)=\mathbf a_h(0)-\bar{\mathbf a}_h$ for the polynomial variables $\mathbf a_h$.

The constant-$\nu$ and constant-$\mu$ tests use 17 and 13 Fourier points, reflection-symmetry reduction, tensor scale $0.5$, and twelve fixed linear substeps. Density is reconstructed classically as $\rho=e^q$ after interpolation to the 129-point reference grid. With the acoustic-mean profile, $E_{\rm fluc}$ falls from $1.04\times10^{-2}$ to $5.81\times10^{-4}$ by order three for constant $\nu$, and from $1.04\times10^{-2}$ to $2.373\times10^{-3}$ by order two for constant $\mu$ (Table~\ref{tab:extension-errors}). The $\mu$ endpoint average improves pooled first-order error to $8.85\times10^{-3}$, but its velocity-component error increases from $1.260\%$ to $1.396\%$. At $\mu,m=2$, the reciprocal constraint $\max_i|s_i e^{q_i}-1|=2.004\times10^{-3}$ and complete-lift LCHS/direct error is $4.793\times10^{-6}$; Fig.~\ref{fig:ns-strong-lchs} shows the reconstructed fields.

\begin{table}[!htbp]
\centering
\TableFont
\caption{Finite-LCHS extension tests using the equidistant Low--Somma rectangle rule. $D_{\rm sym}$ is the active dimension after exact symmetry reduction; $N_Q$ is its exact zero-padded register dimension. $E_{\rm fluc}$ is defined in Eq.~\eqref{eq:extension-fluctuation-error}. All rows use $K_1=K_2=32$ and 385 nodes.}
\label{tab:extension-errors}
\begin{tabular*}{\linewidth}{@{\extracolsep{\fill}}llrrrr@{}}
\toprule
Problem / profile & Grid & $m$ & $D_{\rm sym}$ & $N_Q$ & $E_{\rm fluc}$\\
\midrule
2D Burgers / initial & $9\times9$ & 1 & 13,689 & 16,384 & $3.49\times10^{-5}$\\
2D ZK / initial & $9\times9$ & 1 & 3,564 & 4,096 & $7.03\times10^{-4}$\\
NS, $\nu$ / mean & 17 & 1 & 204 & 256 & $1.04\times10^{-2}$\\
NS, $\nu$ / mean & 17 & 2 & 1,479 & 2,048 & $2.40\times10^{-3}$\\
NS, $\nu$ / mean & 17 & 3 & 9,384 & 16,384 & $5.81\times10^{-4}$\\
NS, $\mu$ / mean & 13 & 1 & 1,810 & 2,048 & $1.04\times10^{-2}$\\
NS, $\mu$ / mean & 13 & 2 & 57,789 & 65,536 & $2.373\times10^{-3}$\\
NS, $\mu$ / endpoints & 13 & 1 & 1,810 & 2,048 & $8.85\times10^{-3}$\\
\bottomrule
\end{tabular*}
\end{table}

\begin{figure}[!htbp]
\centering
\includegraphics[width=\linewidth]{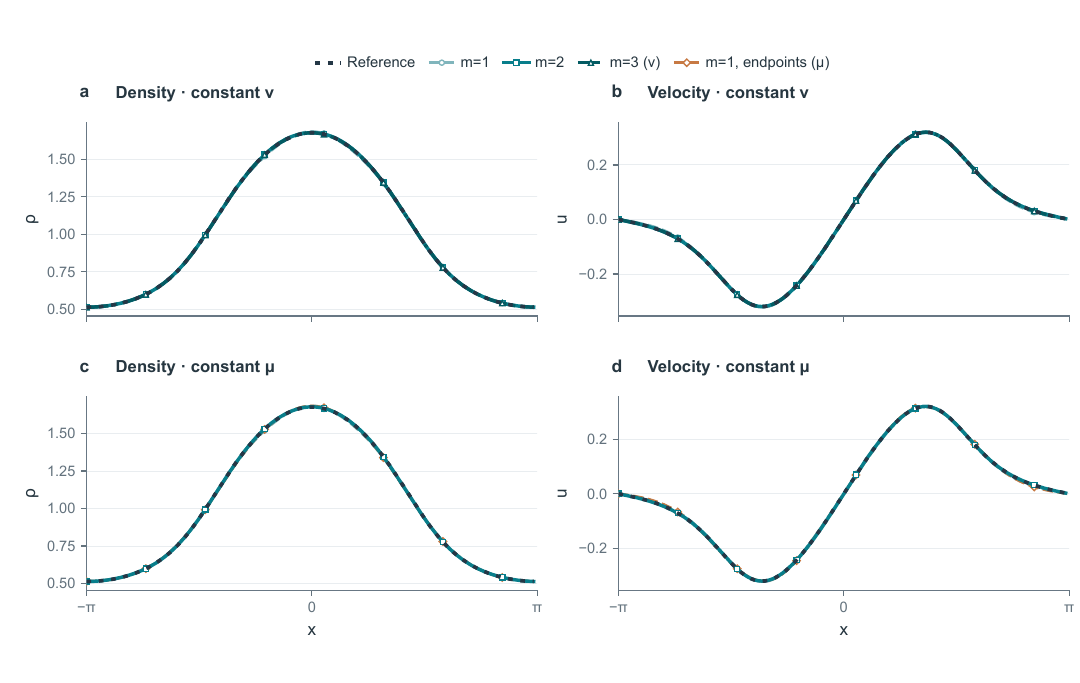}
\caption{One-dimensional strong-compression isothermal Navier--Stokes test for Eq.~\eqref{eq:ns-test} at $T=1.5$. Columns show density and velocity; rows use constant kinematic ($\nu$) and dynamic ($\mu$) viscosity. Dotted black curves are the corresponding 129-point nonlinear references. Colored curves are finite-$K_1/K_2$ Frozen--Oseen HAM-LCHS results at the indicated orders; the $\mu$ row also shows a first-order endpoint-profile result. Frozen profiles use only the linear acoustic predictor.}
\label{fig:ns-strong-lchs}
\end{figure}

\FloatBarrier

\subsection{Two-dimensional isothermal compressible Navier--Stokes equations}
\label{sec:ns2d-tests}

On $[-\pi,\pi)^2$ with periodic boundaries, we take $p=\rho$ and
$\boldsymbol\tau=\nu\rho[\nabla\boldsymbol U+(\nabla\boldsymbol U)^T]$, where $\boldsymbol U=(u,v)^T$ and $\nu$ is the
constant kinematic shear viscosity. The dynamic shear viscosity is
$\nu\rho$ and the second viscosity coefficient is zero. With $q=\log\rho$,
the smooth positive-density equations become the quadratic system
\[
\begin{aligned}
q_t&=-\boldsymbol U\cdot\nabla q-\nabla\cdot\boldsymbol U,\\
\boldsymbol U_t&=-(\boldsymbol U\cdot\nabla)\boldsymbol U-\nabla q\\
&\quad+\nu\{\Delta\boldsymbol U+\nabla(\nabla\cdot\boldsymbol U)
 +[\nabla\boldsymbol U+(\nabla\boldsymbol U)^T]\nabla q\}.
\end{aligned}
\tag{112}\label{eq:ns2d-quadratic}
\]
This retains cross-stress terms; its one-dimensional longitudinal viscosity is $2\nu$, not the coefficient in Section~\ref{sec:ns-tests}.
Both cases start from
\[
q_0=0.18\cos x+0.12\cos y+0.06\cos x\cos y.
\tag{113}
\]
The coupled-compression case uses
$u_0=-0.45\sin x(1+0.2\cos y)$,
$v_0=-0.35\sin y(1+0.15\cos x)$, $(\nu,T)=(0.02,1.2)$.
The vortex--acoustic case uses
$u_0=-0.25\sin x+0.25\sin x\cos y$,
$v_0=-0.20\sin y-0.25\cos x\sin y$, $(\nu,T)=(0.01,1.25)$.

We compare initial-state freezing with the time average of the linear acoustic--viscous predictor
\[
q_t^a=-\nabla\cdot\boldsymbol U^a,\qquad
\boldsymbol U_t^a=-\nabla q^a+\nu\{\Delta\boldsymbol U^a+\nabla(\nabla\cdot\boldsymbol U^a)\},
\tag{114}
\]
initialized at $(q_0,\boldsymbol U_0)$; no nonlinear reference enters the profile. We freeze the full semidiscrete derivative $G=DF_h(\bar{\mathbf a}_h)$ once, with $\mathbf a_h=(q,u,v)$. Reflection symmetry and exact symmetric-product lifting give first-order dimensions 2340 and 4944 on $9\times9$ and $11\times11$ Fourier grids. The runs use tensor scale $0.5$ and twelve fixed linear substeps.

The $31\times31$ nonlinear reference is independently grid-refined. At $9\times9$, acoustic averaging reduces $E_{\rm fluc}$ from $1.974\times10^{-2}$ to $2.522\times10^{-3}$ for compression and from $9.399\times10^{-3}$ to $2.288\times10^{-3}$ for vortex--acoustic interaction. At $11\times11$ the errors are $1.632\times10^{-3}$ and $1.077\times10^{-3}$ (Table~\ref{tab:ns2d-errors}). Figures~\ref{fig:ns2d-compression-fields} and~\ref{fig:ns2d-vortex-fields} show the full fields and signed errors. The plotted fields use finite LCHS; direct propagation is diagnostic only.

\begin{table}[!htbp]
\centering
\TableFont
\caption{First-order finite-LCHS results for two-dimensional NS using 385 equidistant Low--Somma rectangle nodes. $N_Q$ is the padded dimension; mean profiles use the time-averaged linear predictor. $E_{\rm fluc}$ pools the fields: for refined compression the component errors $(\rho,u,v)$ are $(1.29\times10^{-3},1.13\times10^{-2},4.91\times10^{-3})$, so $u$ remains above $1\%$.}
\label{tab:ns2d-errors}
\begin{tabular*}{\linewidth}{@{\extracolsep{\fill}}llrrrr@{}}
\toprule
Case & Profile & Grid & $D_{\rm sym}$ & $N_Q$ & $E_{\rm fluc}$\\
\midrule
Compression & Initial & $9\times9$ & 2,340 & 4,096 & $1.974\times10^{-2}$\\
 & Mean & $9\times9$ & 2,340 & 4,096 & $2.522\times10^{-3}$\\
 & Mean & $11\times11$ & 4,944 & 8,192 & $1.632\times10^{-3}$\\
Vortex--acoustic & Initial & $9\times9$ & 2,340 & 4,096 & $9.399\times10^{-3}$\\
 & Mean & $9\times9$ & 2,340 & 4,096 & $2.288\times10^{-3}$\\
 & Mean & $11\times11$ & 4,944 & 8,192 & $1.077\times10^{-3}$\\
\bottomrule
\end{tabular*}
\end{table}

Across all six runs, complete-lift LCHS/direct discrepancies stay below $4.4\times10^{-8}$. The refined compression case nevertheless has an $u$-component error above $1\%$ (Table~\ref{tab:ns2d-errors}).

\begin{figure}[p]
\centering
\includegraphics[width=\linewidth]{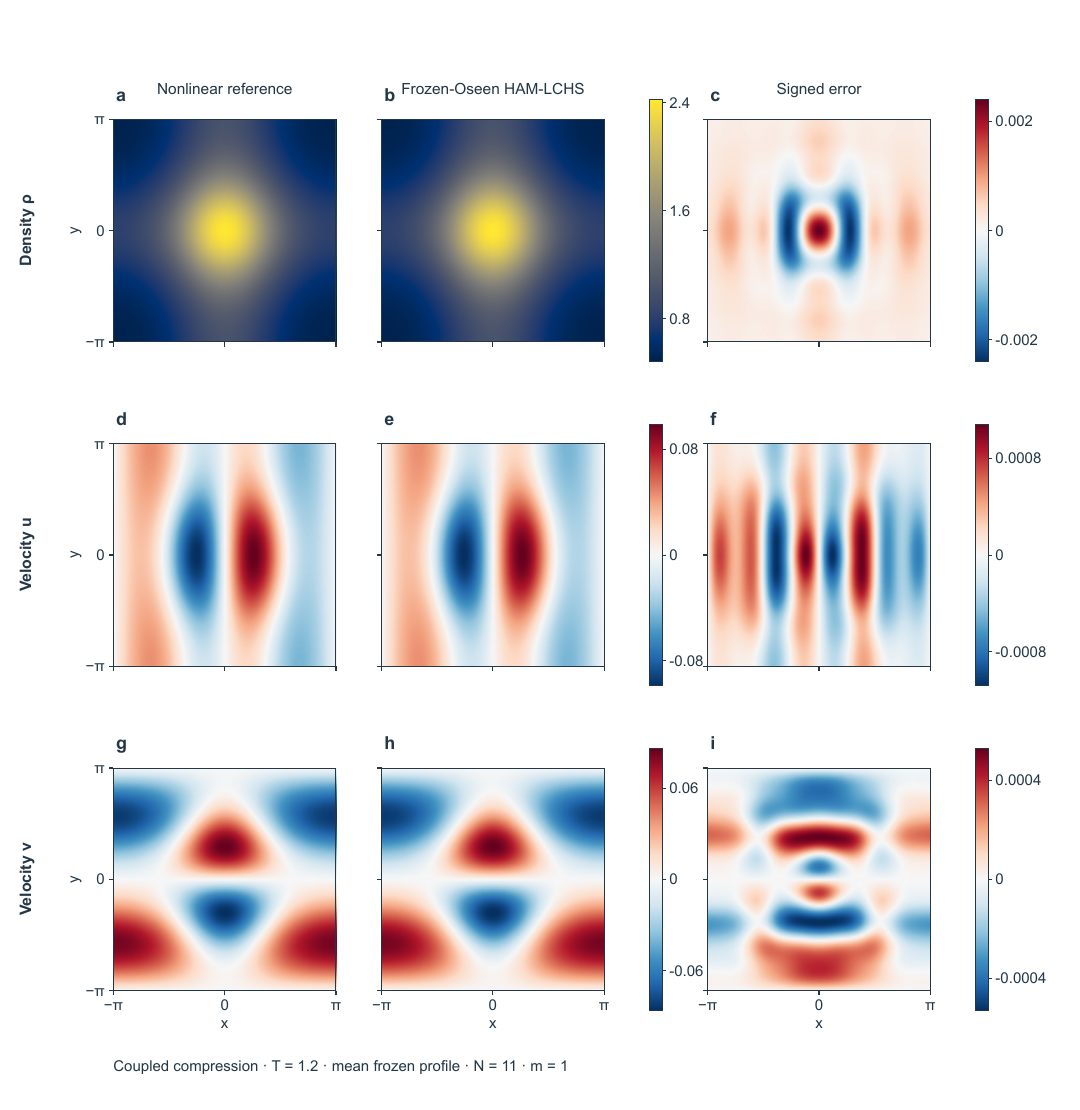}
\caption{Full spatial fields for the coupled-compression case of Eq.~\eqref{eq:ns2d-quadratic} at $T=1.2$. Rows show $\rho$, $u$, and $v$; columns show the $31\times31$ nonlinear reference, the $11\times11$ first-order finite-$K_1/K_2$ Frozen--Oseen HAM-LCHS result with a time-averaged frozen profile, and pointwise signed error $e_\phi=\phi_{\rm LCHS}-\phi_{\rm ref}$ for the displayed field $\phi$. Each reference/result pair shares a scale; errors use separate zero-centered scales. Physical errors are evaluated after Fourier interpolation of $q,u,v$ to the common $31\times31$ grid and reconstruction of $\rho=e^q$; only the separate $128\times128$ resampling is for display.}
\label{fig:ns2d-compression-fields}
\end{figure}

\begin{figure}[p]
\centering
\includegraphics[width=\linewidth]{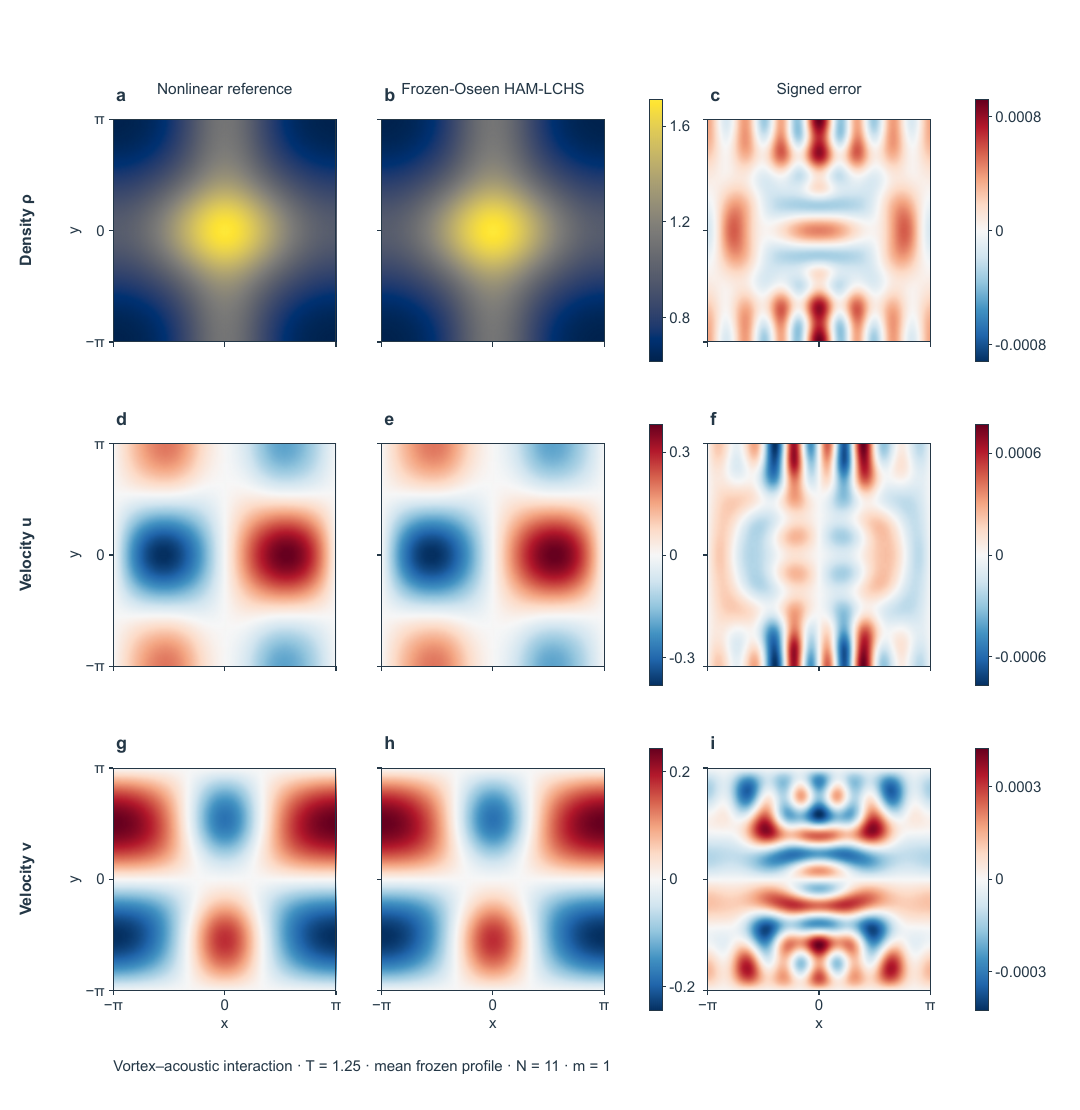}
\caption{Full spatial fields for the vortex--acoustic case of Eq.~\eqref{eq:ns2d-quadratic} at $T=1.25$. Rows show $\rho$, $u$, and $v$; columns show the $31\times31$ nonlinear reference, the $11\times11$ first-order finite-$K_1/K_2$ Frozen--Oseen HAM-LCHS result with a time-averaged frozen profile, and pointwise signed error $e_\phi=\phi_{\rm LCHS}-\phi_{\rm ref}$ for the displayed field $\phi$. Reference/result pairs share a scale; errors use separate zero-centered scales. Physical errors use Fourier interpolation of $q,u,v$ and density reconstruction on the common $31\times31$ grid; separate $128\times128$ resampling is for display only. No fitting to the nonlinear reference or phase alignment is applied.}
\label{fig:ns2d-vortex-fields}
\end{figure}

\FloatBarrier

\subsection{Exploratory two-dimensional incompressible vortex-pair test}
\label{sec:vortex-pair-test}

To test stronger nonlinear transport beyond the smooth manufactured cases, we consider a pair of like-signed vortices on $[-\pi,\pi)^2$~\cite{melander1988symmetric}. In vorticity--streamfunction variables,
\[
\omega_t+\boldsymbol V(\omega)\cdot\nabla\omega=\nu\Delta\omega,\qquad
-\Delta\psi=\omega,\qquad \boldsymbol V(\omega)=(\psi_y,-\psi_x)^T,
\tag{115}\label{eq:vortex-pair}
\]
with zero spatial mean and periodic boundaries. We set $\nu=0.02$, $T=3$, and $\omega_0=g-\langle g\rangle$, where
\[
g(x,y)=1.5\sum_{\sigma=\pm1}
\exp\!\left\{2[\cos(x-0.8\sigma)+\cos y-2]\right\}.
\tag{116}
\]
\begin{figure}[p]
\centering
\includegraphics[width=\linewidth]{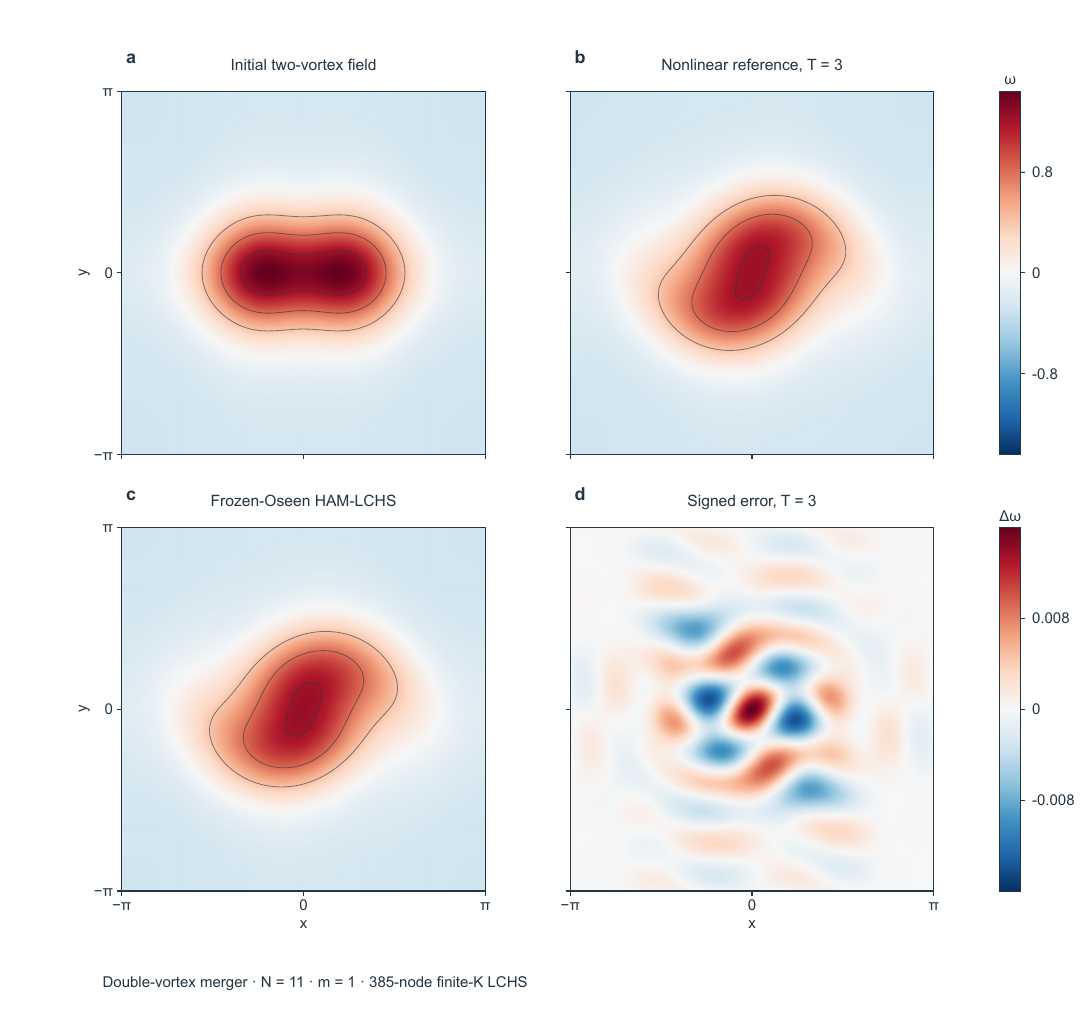}
\caption{Exploratory vortex-pair test for Eq.~\eqref{eq:vortex-pair} at $T=3$. Panels show (a) initial vorticity, (b) the $31\times31$ nonlinear reference, (c) the $11\times11$ first-order finite-$K_1/K_2$ Frozen--Oseen HAM-LCHS result, and (d) pointwise signed error $e_\omega=\omega_{\rm LCHS}-\omega_{\rm ref}$. Panels (a)--(c) share a color scale; (d) uses a separate zero-centered scale. Periodic resampling is for display only.}
\label{fig:ns2d-vortex-pair}
\end{figure}

We freeze $\bar\omega=\omega_0$, giving $\mathcal Gz=\nu\Delta z-\boldsymbol V(\omega_0)\cdot\nabla z-\boldsymbol V(z)\cdot\nabla\omega_0$. Rotational symmetry and zero mean reduce the $11\times11$ system to 60 physical coordinates and a 2010-dimensional first-order lift, padded to 2048. With twelve fixed linear substeps and tensor scale $0.1$, finite LCHS gives final vorticity error $7.33\times10^{-3}$ against the $31\times31$ nonlinear reference; the $11\times11$ spatial discrepancy is $5.95\times10^{-3}$ and complete-lift LCHS/direct discrepancy $1.13\times10^{-7}$. Figure~\ref{fig:ns2d-vortex-pair} shows the two peaks merging; this remains an exploratory first-order test.

\FloatBarrier

\subsection{Matched linear-observable comparison}\label{matched-observable-comparison}

We compare two fixed Fourier-mode outputs of each saved finite-LCHS field with direct nonlinear integration on the same semidiscrete system. For a row-orthonormal observable matrix $W$ acting on the computational variables $\mathbf a_h$, the relative error is $E_{\rm obs}=\|W\mathbf a_h-W\mathbf a_{h,\rm ref}\|_2/\|W\mathbf a_{h,\rm ref}\|_2$. Burgers uses two sine modes; constant-$\mu$ NS uses one $q=\log\rho$ cosine mode and one velocity sine mode.

\begin{table}[!htbp]
\centering
\caption{Same-output comparison at the existing Burgers \(n=5,m=5,T=1\) and constant-\(\mu\) NS \(N_x=13,m=2,T=1.5\) settings. DOP853 uses relative tolerance \(10^{-3}\). Classical times include model construction, integration and observable evaluation; they are medians of five warmed runs with one BLAS thread, not quantum runtimes. Reference generation and prior profile calibration are excluded.}
\label{tab:matched-observable}
\TableFont
\begin{tabular*}{\linewidth}{@{\extracolsep{\fill}}lrrr@{}}
\toprule
Case & \(E_{\rm obs}\), finite LCHS & \(E_{\rm obs}\), DOP853 & Time (ms)\\
\midrule
Burgers & \(4.985\times10^{-4}\) & \(3.517\times10^{-6}\) & 0.884\\
1D NS & \(3.582\times10^{-5}\) & \(5.384\times10^{-7}\) & 3.401\\
\bottomrule
\end{tabular*}
\end{table}

At these small-grid tolerances, DOP853 is more accurate and faster than explicit lifted propagation (Table~\ref{tab:matched-observable}). The comparison does not establish a quantum speedup or make explicit-lift storage a necessary classical baseline.
\FloatBarrier
\endgroup

\section{Conclusions and discussion}\label{conclusions-and-discussion}
The central FOQHAM framework connects transport-aware auxiliary-operator selection to the finite linear representation of nonlinear PDE dynamics. Freezing the full Fr\'echet derivative retains transport and profile-gradient coupling, with an exact residual and a locally quadratic remainder. Product lifting then closes each prescribed finite-order polynomial hierarchy as a fixed affine linear system. LCHS propagates the linear system without outer homotopy iterations and profile updates. Local approximation and profile-sensitivity bounds characterize conditions under which the chosen auxiliary can improve the approximation.

Numerical experiments on Burgers, Korteweg–de Vries, Zakharov–Kuznetsov and isothermal compressible Navier–Stokes equations support accurate non-iterative approximations in the tested regimes. Matched auxiliary controls connect a reduced homotopy order to a smaller lifted representation in a Burgers case. The KdV comparison instead favors mean transport over the full frozen derivative at the prescribed accuracy. These results identify auxiliary selection as a joint approximation and resource-design problem, involving the nonlinear remainder, propagator growth and representation size.

The framework complements quantum linear-evolution algorithms by organizing nonlinear approximation before propagation. Under the stated matrix-access, normalization, signal and readout conditions, its resource analysis provides a basis for quantum implementation. Structured compressed encodings and independently validated profiles are concrete next steps towards larger, spatially resolved systems. More broadly, this work shows how physically informed auxiliary operators can connect nonlinear PDE approximation to quantum linear evolution, providing a systematic foundation for further quantum algorithm development.

\section*{Acknowledgments}
Human authors came up with the Frozen-Oseen operator. OpenAI GPT-5.6 and 6 assisted with language polishing in chapter 2-4, theorem proofs in chapter 2, mathematical and citation checks, numerical tests' code improving. The final proofs, algorithms, numerical tests' code was independently written , reviewed, verified and approved by the human authors, who take full responsibility for the work.

\bibliographystyle{elsarticle-num}
\bibliography{oseen_qham_jcp_references}

\end{document}